%% file: main.tex
\documentclass[hidelinks,onefignum,onetabnum]{siamart220329}

\usepackage{amsfonts}

\usepackage[inline]{enumitem}
\usepackage{macros}

\usepackage{crossreftools}
\makeatletter
\def\bstctlcite{\@ifnextchar[{\@bstctlcite}{\@bstctlcite[@auxout]}}
\def\@bstctlcite[#1]#2{\@bsphack
  \@for\@citeb:=#2\do{%
    \edef\@citeb{\expandafter\@firstofone\@citeb}%
    \if@filesw\immediate\write\csname #1\endcsname{\string\citation{\@citeb}}\fi}%
  \@esphack}
\makeatother 

\usepackage[T2A]{fontenc}

\usepackage{hyphenat}
\newsiamremark{remark}{Remark}
\newsiamremark{hypothesis}{Hypothesis}
\crefname{hypothesis}{Hypothesis}{Hypotheses}
\newsiamthm{claim}{Claim}
\newsiamthm{example}{Example}

\headers{Distributional \texorpdfstring{\lowercase{$k$}}{k}-Plane Transform}{R. Parhi and M. Unser}

\title{Distributional Extension and Invertibility of the \texorpdfstring{\lowercase{$k$}}{k}-Plane Transform and Its Dual\thanks{\vspace{-1em}
\funding{This work was supported by the Swiss National Science Foundation under Grant 200020\_219356 / 1.}}}

\author{Rahul Parhi\thanks{Biomedical Imaging Group, \'Ecole polytechnique f\'ed\'erale de Lausanne, CH-1015 Lausanne, Switzerland (\email{rahul.parhi@epfl.ch}, \email{michael.unser@epfl.ch}).}
\and Michael Unser\footnotemark[2]}

\ifpdf
\hypersetup{
  pdftitle={Distributional \texorpdfstring{\lowercase{$k$}}{k}-Plane Transform},
  pdfauthor={R. Parhi and M. Unser}
}
\fi

\begin{document}

\bstctlcite{IEEEexample:BSTcontrol}

\maketitle

\begin{abstract}
    \input{sections/abstract}
\end{abstract}



\input{sections/introduction}
\input{sections/preliminaries}
\input{sections/k-plane}

\input{sections/k-plane-Schwartz}
\input{sections/distributional}
\input{sections/Fourier-slice}
\input{sections/compatible}

\input{sections/specific-transforms}
\input{sections/inverse-problems}

\appendix

\input{appendix/iso-projector}

\section*{Acknowledgment}
The authors would like to thank the anonymous referee for their careful reading of the manuscript and for providing concrete suggestions to improve the clarity of the exposition. Their feedback led to substantial improvement in the content and presentation of this work.

\bibliographystyle{IEEEtranS}
\bibliography{ref}

\end{document}

%% file: sections/abstract.tex
We investigate the distributional extension of the $k$-plane transform in $\R^d$ and of related operators. We parameterize the $k$-plane domain as the Cartesian product of the Stiefel manifold of orthonormal $k$-frames in $\R^d$ with $\R^{d-k}$. This parameterization imposes an isotropy condition on the range of the $k$-plane transform which is analogous to the even condition on the range of the Radon transform. We use our distributional formalism to investigate the invertibility of the dual $k$-plane transform (the ``backprojection'' operator). We provide a systematic construction (via a completion process) to identify Banach spaces in which the backprojection operator is invertible and present some prototypical examples. These include the space of isotropic finite Radon measures and isotropic $L^p$-functions for $1 < p < \infty$. Finally, we apply our results to study a new form of regularization for inverse problems.

%% file: sections/introduction.tex
\section{Introduction}
The $k$-plane transform of a function $f: \R^d \to \R$, denoted by $\RadonOp_k\curly{f}$, is defined on the space of $k$-planes (affine subspaces of dimension $k$) in $\R^d$ for $1 \leq k < d$. The value of $\RadonOp_k\curly{f}$ at a given $k$-plane is the integral of $f$ over that $k$-plane. The two extremes of this transform correspond to the Radon transform ($k = (d - 1)$) and the X-ray transform ($k = 1$). Intermediate dimensions are also of interest such as $k = (d-2)$, which corresponds to the measurement model of nuclear magnetic resonance imaging~\cite{SmithRadiographs,MR0548770}. These transforms have had considerable interest in both pure and applied mathematics following the seminal work of Radon~\cite{Radon} and John~\cite{JohnBestimmung,JohnAbhangigkeiten}. The Radon transform, in particular, has had profound impacts in the study of partial differential equations (PDEs)~\cite{JohnBook}, computerized tomography~\cite{LoganCT}, the theory of wavelets, ridgelets, and shearlets~\cite{BartolucciRadonShearlet,CandesPhD,CandesRidgelets,DonohoRidgelets,KostadinovaRidgelets,kostadinova2015ridgelet,RubinRidgelets,shenouda2023continuous,SonodaRidgelets}, and, more recently, neural networks~\cite{BartolucciRKBS,KurkovaEstimates,OngieRadon,ParhiShallowRepresenter,ParhiDeepRepresenter,ParhiSPM,ParhiMinimax}. Moreover, there is a large body of work that has investigated properties of the Radon transform. These include characterizations of its range on various spaces of functions~\cite{HelgasonDuality,HelgasonRadonEuclidean,JohnBestimmung,LudwigRadon}, distributional extensions~\cite{KostadinovaRidgelets,LudwigRadon}, the invertibility of the dual transform (the ``backprojection'' operator)~\cite{KostadinovaRidgelets,RammBackprojection,SolmonDualRadon}, and several others~\cite{AlbertiRadonDualPairs,KatsevichResolution,MarkoeInvertibility}.

The $k$-plane transform, sometimes referred to as the $k$-dimensional Radon transform or the Radon--John transform, is the most obvious generalization of the Radon transform. While there are works devoted to the $k$-plane transform that include inversion formulas~\cite{Fuglede,KeinertInversion,RubinInversion,RubinKPlaneInversion,RubinInversionFormulas}, range characterizations~\cite{GonzalezRange,GrinbergRangeKPlane,MarkoeAnalyticTomo,petrov-en,RichterRangeKPlane,PetrovCavalieri1,PetrovCavalieri2,PetrovCavalieri3}, as well as its connections to wavelets~\cite{DonohoKPlaneRidgelets,RubinInversion,RubinConvolutionBackprojection}, far less is known about the $k$-plane transform than the Radon transform.

This paper brings the following new results:
\begin{enumerate}
    \item The extension of the $k$-plane transform and of related operators to distributions (\cref{sec:dist-ext}).
    \item The formulation of an extended version of the Fourier slice theorem that is applicable to any tempered distribution (\cref{thm:Fourier-slice}), a result that is new, to the best of our knowledge, for all $1 \leq k < d$ (which includes the Radon transform). \Cref{thm:Fourier-slice} also yields a general form of the intertwining properties for the $k$-plane transform (\cref{thm:intertwine}).
    \item The systematic construction (via a completion process) of Banach spaces over which the dual $k$-plane transform (backprojection operator) is guaranteed to be invertible (\cref{thm:k-plane-systematic}).
    \item The identification of prototypical examples of such $k$-plane-domain Banach spaces. These examples recover and extend many known results regarding the invertibility of the backprojection operator (see \cref{rem:dual-Radon} for a detailed discussion). They include the space of isotropic finite Radon measures and isotropic $L^p$-functions for $1 < p < \infty$
    (\cref{cor:backprojection-invertible}).
\end{enumerate}
The primary technical contribution of this paper is the careful treatment of the extension of the $k$-plane transform and related operators to distributions. The main difficulty is that the range of the $k$-plane transform on Schwartz functions is not the full space of Schwartz functions on the $k$-plane domain. It is instead a subspace of Schwartz functions over which certain consistency conditions are met. As it turns out, when we extend the $k$-plane transform by duality, the resulting distribution spaces are actually \emph{equivalence classes} of distributions. This introduces an inherent ambiguity and prevents one to identify the $k$-plane transform of a distribution uniquely. To deal with this non-uniqueness issue, we construct $k$-plane-domain Banach spaces over which the backprojection operator is invertible. This then yields a setting where the distributional (filtered) $k$-plane transform can be \emph{uniquely determined} (\cref{thm:k-plane-systematic}). 

\section{Motivation}
    The primary motivation for the proposed investigation is to study \emph{$k$-plane ridges}. These are functions that map $\R^d \to \R$ and that are constant along $k$-planes with a fixed orientation. They can be written as
    \begin{equation}
        \vec{x} \mapsto r(\mat{A}\vec{x}), \quad \vec{x} \in \R^d,
        \label{eq:k-plane-ridge-intro}
    \end{equation}
    where $r: \R^{d-k} \to \R$ is the \emph{profile} and $\mat{A} \in \R^{(d-k) \times d}$ is the \emph{orientation} of the $k$-plane ridge.
    Observe that \cref{eq:k-plane-ridge-intro} is constant along the $k$-planes $\mat{A}\vec{x} = \vec{c}$, where $\vec{c} \in \R^{d-k}$, whenever $\mat{A}$ is full rank.
    
    Early works that have used the $k$-plane transform to study $k$-plane ridges go back to the 1990s and have been spearheaded by  Donoho~\cite{DonohoKPlaneRidgelets} and Rubin~\cite{RubinConvolutionBackprojection}. There, the authors propose wavelet-like systems that interpolate between wavelets and ridgelets, as $k$ varies. These systems, referred to as $k$-plane ridgelets, are built from $k$-plane ridges with $(d-k)$-variate wavelet profiles. Crucially, they rely on wavelet analysis in the $k$-plane domain. Donoho's motivation was to propose a system that could capture and exploit low-dimensional structures found in high-dimensional functions and data, with the purpose of understanding the curse of dimensionality in high-dimensional data-analysis problems~\cite{DonohoCurse}.

    More recently, $k$-plane ridges, sometimes referred to as \emph{multi-index models} or \emph{$(d-k)$-sparse functions}, have been studied in the context of dimensionality reduction~\cite{dalalyan2008new,fukumizu2004dimensionality} and as a tool to understand why neural networks work well in high-dimensional problems---whereas kernel methods do not~\cite{Bach,ghorbani2020neural}. In particular, the distributional filtered Radon transform ($k = (d-1)$) and ridgelet analysis have been used to characterize a quantitative gap between the performance of neural networks and kernel methods via nonlinear and linear minimax rates for nonparametric function estimation~\cite{ParhiMinimax}. Notably, the authors of~\cite{ParhiMinimax} prove that, for certain classes of problems, neural networks enjoy estimation error rates that do not grow with the input dimension, whereas the estimation error rates for kernel methods \emph{necessarily} suffer the curse of dimensionality.

    To summarize, we see that $k$-plane ridges and the $k$-plane transform are powerful tools to understand and quantify the curse of dimensionality. Therefore, we believe that the results of this paper provide the mathematical framework to further the investigation of this hallmark. Indeed, in \cref{sec:specific-transforms} we use our formalism to prove that the application of the distributional filtered $k$-plane transform to a $k$-plane ridge extracts the profile function and orientation matrix (see \cref{thm:k-plane-ridges}). This property allows for one to analyze combinations of $k$-plane ridges and to access the $(d-k)$-variate profiles directly, without knowledge of the underlying orientation matrices. Said differently, the filtered $k$-plane transform automatically captures intrinsic low-dimensional structures in high-dimensional functions. We exemplify the utility of this result in \cref{sec:inverse-problems}, where we investigate the resolution of inverse problems with new forms of $k$-plane-domain regularization.

%% file: sections/preliminaries.tex
\section{Preliminaries and Notation}
The Schwartz space of smooth and rapidly decreasing functions on $\R^d$ is denoted by $\Sch(\R^d)$. Its continuous dual is the space $\Sch'(\R^d)$ of tempered distributions. We let $L^p(\R^d)$ denote the Lebesgue space for $1 \leq p \leq \infty$ and $L^1_\temp(\R^d) = L^1_\loc(\R^d) \cap \Sch'(\R^d)$ the space of tempered locally integrable functions. The Banach space of continuous functions vanishing at $\pm \infty$ on $\R^d$ equipped with the $L^\infty$-norm is denoted by $C_0(\R^d)$. By the Riesz--Markov--Kakutani representation theorem~\cite[Chapter~7]{FollandRA}, its continuous dual can be identified with the Banach space of finite Radon measures, denoted by $\M(\R^d)$. Since $\Sch(\R^d) \dembed C_0(\R^d)$, where $\dembed$ denotes a dense and continuous embedding, we have by duality that $\M(\R^d) \cembed \Sch'(\R^d)$, where $\cembed$ denotes a continuous embedding. In the case of a compact domain $\Omega$, we have $\M(\Omega) = \paren{C(\Omega)}'$, where $C(\Omega)$ denotes the space of continuous functions on $\Omega$.

The Fourier transform of $\varphi \in \Sch(\R^d)$ is defined as
\begin{equation}
    \hat{\varphi}(\vec{\xi}) \coloneqq \FourierOp_d\curly{\varphi}(\vec{\xi}) = \int_{\R^d} \varphi(\vec{x}) e^{-\imag \vec{\xi}^\T\vec{x}} \dd\vec{x}, \quad \vec{\xi} \in \R^d,
\end{equation}
where $\imag^2 = -1$. Consequently, the inverse Fourier transform of $\hat{\varphi} \in \Sch(\R^d)$ is given by
\begin{equation}
    \FourierOp_d^{-1}\curly{\hat{\varphi}}(\vec{x}) = \frac{1}{(2\pi)^d} \int_{\R^d} \hat{\varphi}(\vec{\xi}) e^{\imag \vec{\xi}^\T\vec{x}} \dd\vec{\xi}, \quad \vec{x} \in \R^d.
\end{equation}
These operators are extended to act on $\Sch'(\R^d)$ by duality. Any continuous linear shift-invariant (LSI) operator $\LOp: \Sch(\R^d) \to \Sch'(\R^d)$ is a convolution operator specified by a unique kernel $h \in \Sch'(\R^d)$ such that $\LOp \varphi = h * \varphi$. Such operators can also be specified in the Fourier domain by
\begin{equation}
    \LOp \varphi = \FourierOp_d^{-1}\curly{\hat{L} \hat{\varphi}},
\end{equation}
where $\hat{L} \in \Sch'(\R^d)$ is the Fourier transform of the kernel $h \in \Sch'(\R^d)$. The tempered distribution $h$ is the \emph{impulse response} of $\LOp$ and the tempered distribution $\hat{L}$ is the Fourier symbol or \emph{frequency response} of $\LOp$. We shall generally use upright, roman letters for LSI operators and use the italic variant with a hat to denote its frequency response. Finally, given a space $\X$ and a norm $\norm{\dummy}$, the completion of $\X$ in $\norm{\dummy}$ is a Banach space, denoted by $\cl{(\X, \norm{\dummy})}$.

%% file: sections/k-plane.tex
\section{Parameterization of the \texorpdfstring{$k$}{k}-Plane Transform}
The Radon transform is a special case of the $k$-plane transform.
It corresponds to integrals over affine subspaces of dimension
$(d-1)$ in $\R^d$ (i.e., hyperplanes).
The classical integral formula of the Radon transform of $f \in L^1(\R^d)$ is
\begin{equation}
  \Radon{f}(\vec{\alpha}, t) = \int_{\vec{\alpha}^\perp} f(\vec{x} + \vec{\alpha}t) \dd\vec{x},  \quad (\vec{\alpha}, t) \in \cyl,
  \label{eq:Radon-perp}
\end{equation}
where $\vec{\alpha}^\perp \coloneqq \curly{\vec{x} \in \R^d \st \vec{\alpha}^\T\vec{x} =
0}$ and $\Sph^{d-1} = \curly{\vec{x} \in \R^d \st \norm{\vec{x}}_2 = 1}$ is the unit sphere in $\R^d$. This transform integrates $f$ along the hyperplane $\curly{\vec{x} \in \R^d \st \vec{\alpha}^\T\vec{x} = t}$. In the case of $\varphi \in \Sch(\R^d)$, this transform is equivalently specified by
\begin{equation}
  \Radon{\varphi}(\vec{\alpha}, t) = \int_{\R^d} \varphi(\vec{x})
  \delta(\vec{\alpha}^\T\vec{x} - t) \dd\vec{x},
  \label{eq:Radon-dirac}
\end{equation}
where $\delta$ is the univariate Dirac impulse and the integral is understood as the action of $\delta(\vec{\alpha}^\T(\dummy) - t) \in \Sch'(\R^d)$ on $\varphi \in \Sch(\R^d)$. Since $(\vec{\alpha}, t)$ and $(-\vec{\alpha}, -t)$ specify the same hyperplane, we see that the Radon transform imposes an even symmetry.

Many authors work with the $k$-plane transform in $\R^d$ somewhat abstractly by
considering integrals of a function $f: \R^d \to \R$ over elements of the
affine Grassmann manifold which
is precisely the manifold of $k$-planes in $\R^d$ (see~\cite{MarkoeAnalyticTomo} and the references therein). This abstract formulation
does not recover the usual specifications of the Radon transform given
by \cref{eq:Radon-perp,eq:Radon-dirac} in an immediate way.
Following the work of Petrov~\cite{PetrovCavalieri1,PetrovCavalieri2,PetrovCavalieri3}, we can alternatively work with a concrete parameterization of the affine Grassman manifold with matrices and vectors. Here, we specify a $k$-plane in $\R^d$ by the set of all $\vec{x} \in \R^d$ that satisfy the system of equations
\begin{equation}
  \mat{A}\vec{x} = \vec{t},
\end{equation}
where $\mat{A} \in \R^{(d-k) \times d}$, and $\vec{t} \in \R^{d-k}$. Let us write $\mat{A}$ row-wise as
\begin{equation}
  \mat{A} = \begin{bmatrix}\vec{\alpha}_1^\T \\ \rowsvdots \\ \vec{\alpha}_{d-k}^\T \end{bmatrix}
\end{equation}
While Petrov only requires that the rows of $\mat{A}$ are linearly independent, we additionally impose, without loss of generality, that the rows are orthonormal, with 
$\vec{\alpha}_n \in \Sph^{d-1}$, $n = 1, \ldots, (d-k)$, and
$\vec{\alpha}_m^\T\vec{\alpha}_n = \delta[m - n]$, where $\delta[\dummy]$
denotes the Kronecker impulse. Accordingly, $\mat{A}$
is an element of the Stiefel manifold of orthonormal $k$-frames
in $\R^d$ specified by
\begin{equation}
  V_{d-k}(\R^d) = \curly{\mat{A} \in \R^{(d-k) \times d} \st \mat{A} \mat{A}^\T
  = \mat{I}_{d-k}}.
  \label{eq:Stiefel-manifold}
\end{equation}
While we do not explicitly require the manifold structure of
$V_{d-k}(\R^d)$ in this paper, we do require the ability to integrate on $V_{d-k}(\R^d)$ with
respect to its Haar measure, which we denote by $\sigma$.

With this parameterization, we specify the $k$-plane transform of $f \in L^1(\R^d)$ as
\begin{equation}
  \RadonOp_k\curly{f}(\mat{A}, \vec{t}) = \int_{\mat{A}^\perp} f(\vec{x} +
  \mat{A}^\T\vec{t})
  \dd\vec{x},  \quad (\mat{A}, \vec{t}) \in
  (V_{d-k}(\R^d), \R^{d-k}),
  \label{eq:k-plane-perp}
\end{equation}
where $\mat{A}^\perp \coloneqq \curly{\vec{x} \in \R^d \st \mat{A}\vec{x} =
\vec{0}}$. This transform integrates $f$ along the $k$-plane $\curly{\vec{x} \in \R^d \st \mat{A}\vec{x} = \vec{t}}$. In the case of $\varphi \in \Sch(\R^d)$, this transform is equivalently specified by
\begin{equation}
  \RadonOp_k\curly{f}(\mat{A}, \vec{t}) = \int_{\R^d} \varphi(\vec{x})
  \delta(\mat{A}\vec{x} - \vec{t}) \dd\vec{x},
  \label{eq:k-plane-Dirac}
\end{equation}
where $\delta$ is the $(d-k)$-variate Dirac impulse and the integral is understood as the action of $\delta(\mat{A}(\dummy) - \vec{t}) \in \Sch'(\R^d)$ on $\varphi \in \Sch(\R^d)$. In \cref{defn:k-plane-ridge-distribution,rem:k-plane-Dirac}, we expose the precise meaning of the tempered distribution $\delta(\mat{A}(\dummy) - \vec{t})$. We see that \cref{eq:k-plane-perp,eq:k-plane-Dirac} immediately recover the specifications of the Radon transform in \cref{eq:Radon-perp,eq:Radon-dirac} when $k = (d - 1)$.  Further, the dual transform (often called the ``backprojection'') of $g \in L^\infty(V_{d-k}(\R^d) \times \R^{d-k})$ is given by
\begin{equation}
  \RadonOp_k^*\curly{g}(\vec{x}) = \int_{V_{d-k}(\R^d)} g(\mat{A},
  \mat{A}\vec{x}) \dd\sigma(\mat{A}), \quad \vec{x} \in \R^d.
\end{equation}
Furthermore, since we impose that the rows of $\mat{A}$ are orthonormal, we have that $(\mat{A}, \vec{t})$ and $(\mat{U}\mat{A}, \mat{U}\vec{t})$ define the same $k$-plane, for any orthogonal transformation $\mat{U} \in \mathrm{O}_{d-k}(\R)$ (the orthogonal group in dimension $(d-k)$). The main advantage of the proposed parameterization is that it will allow us to identify the symmetries of $k$-plane domain as ``isotropic'' symmetries.

%% file: sections/k-plane-Schwartz.tex
\section{The \texorpdfstring{$k$}{k}-Plane Transform of Schwartz Functions}
Let $\Xi_k = V_{d-k}(\R^d) \times \R^{d-k}$ denote the $k$-plane domain. We exclusively use the variables $(\mat{A}, \vec{t}) \in \Xi_k$ to index the $k$-plane domain. The space of Schwartz functions on $\Xi_k$, denoted by $\Sch(\Xi_k)$, is the space of smooth functions that are rapidly decreasing in the $\vec{t} \in \R^{d-k}$ variable~\cite{GonzalezRange}. More specifically, we have that $\Sch(\Xi_k) = C^\infty(V_{d-k}(\R^d)) \,\hat{\otimes}\, \Sch(\R^{d-k})$, where $\hat{\otimes}$ denotes the topological tensor product, which is the completion of the algebraic tensor product with respect to the projective topology~\cite[Chapter~43]{TrevesTVS}. We have the following classical result regarding the continuity and invertibility of the $k$-plane transform.
\begin{theorem}[{see~\cite{GelfandIntegralGeometry,GonzalezRange,KeinertInversion,RubinInversion,SmithRadiographs,SolmonXRay}}] \label{thm:cont-inv}
  The operator $\RadonOp_k$ continuously maps $\Sch(\R^d)$ into $\Sch(\Xi_k)$.
  Moreover,
  \begin{equation}
    \RadonOp_k^* \KOp_{d-k} \RadonOp_k = c_{d,k} (-\Delta_d)^\frac{k}{2} \RadonOp_k^*
    \RadonOp_k = c_{d,k} \RadonOp_k^* \RadonOp_k (-\Delta_d)^{\frac{k}{2}} = \Id
  \end{equation}
  on $\Sch(\R^d)$, with
  \begin{equation}
    c_{d,k} = \frac{1}{(2\pi)^k} \frac{\abs{\Sph^{k-1}}}{\abs{\Sph^{d-k-1}}} \frac{1}{\sigma(V_{d-k-1}(\R^{d-1}))},
    \label{eq:k-plane-constant}
  \end{equation}
  where $\abs{\dummy}$ denotes the surface area and $\sigma(\dummy)$ is the Haar measure of the Stiefel manifold.
  The underlying operators are the $d$-variate Laplacian operator $\Delta_d$ and the filtering operator\footnote{In computed tomography (CT), this filter is referred to as the backprojection filter found in the filtered backprojection algorithm for CT image reconstruction.} $\KOp_{d-k} = c_{d,k} (-\Delta_{d-k})^{k/2}$, where $\Delta_{d-k}$ denotes the $(d-k)$-variate Laplacian applied to the $\vec{t} \in \R^{d-k}$ variable. The filtering operator is equivalently specified by the frequency response $\hat{K}_{d-k}(\vec{\omega}) = c_{d,k} \norm{\vec{\omega}}_2^k$, $\vec{\omega} \in \R^{d-k}$.
\end{theorem}
\begin{remark}
    We shall often use subscripts on operators to remind the reader of the number of variables of the functions or distributions the operator applies to.
\end{remark}

Using the expression for the volume of the Stiefel manifold from~\cite[p. 9]{EdelmanVolume}, we simplify the constant in \cref{eq:k-plane-constant} to
\begin{equation}
    c_{d,k} = \frac{1}{(2\pi)^k} \frac{\abs{\Sph^{k-1}}}{\abs{\Sph^{d-k-1}}} \frac{1}{\prod_{n=k}^{d-1} \abs{\Sph^{n-1}}}.
\end{equation}
We then strengthen the result of \cref{thm:cont-inv} by invoking an extended version of the open mapping theorem, where we use the key property
that the range of the $k$-plane transform $\Sch_k =
\RadonOp_k\paren*{\Sch(\R^d)}$ is a closed subspace of the Fr\'echet space
$\Sch(\Xi_k)$~\cite[Chapter~4]{MarkoeAnalyticTomo}.
\begin{corollary} \label[corollary]{cor:homeo}
  The operator $\RadonOp_k: \Sch(\R^d) \to \Sch_k$ is a homeomorphism with
  inverse $\RadonOp_k^{-1} = \RadonOp_k^*\KOp_{d-k}: \Sch_k \to
  \Sch(\R^d)$.
\end{corollary}
\begin{proof}
  From \cref{thm:cont-inv}, we know that $\RadonOp_k: \Sch(\R^d) \to \Sch_k$ is
  a continuous bijection with inverse given by $\RadonOp_k^{-1} =
  \RadonOp_k^*\KOp_{d-k}$. Since $\Sch(\R^d)$ and $\Sch_k$ are Fr\'echet
  spaces, where $\Sch_k$ inherits its Fr\'echet topology from 
  $\Sch(\Xi_k)$, by the
  open mapping theorem for Fr\'echet spaces~\cite[Theorem~2.11]{RudinFA}, we
  know that $\RadonOp_k$ is an open map.  Therefore, for $\RadonOp_k^*\KOp_{d-k}$,
  the inverse images of open sets are open, which implies the continuity of $\RadonOp_k^*\KOp_{d-k}$.
\end{proof}

Let $\Sch_H$ denote the subspace of Schwartz functions on $\Xi_k$ that satisfy the so-called Helgason consistency conditions (see~\cite[Chapter~4]{MarkoeAnalyticTomo} for a precise definition). In the case of the Radon transform, these consistency conditions exactly characterize its range ($\Sch_H = \Sch_{d-1}$). But, when $k < (d - 1)$, it turns out that $\Sch_H \neq \Sch_k$~\cite{GonzalezRange}. In fact, $\Sch_k$ is a subspace of $\Sch_H$ whose members satisfy a system of PDEs~\cite{GonzalezRange,MarkoeAnalyticTomo}. By considering a different set of consistency conditions, Petrov could show in~\cite{PetrovCavalieri1,PetrovCavalieri2} (see~\cite{petrov-en} for an English translation of~\cite{PetrovCavalieri1}) that the range has an equivalent characterization which does not explicitly involve PDEs~\cite[Remark~6]{PetrovCavalieri1}. This is summarized in \cref{thm:k-plane-range}.

\begin{theorem}[\cite{petrov-en,PetrovCavalieri1,PetrovCavalieri2}] \label{thm:k-plane-range}
  A function $\psi$ is the $k$-plane transform of a function $\varphi \in
  \Sch(\R^d)$ (i.e., $\psi \in \Sch_k$) if and only if the following conditions hold:
  \begin{enumerate}
    \item $\psi \in \Sch(\Xi_k)$;
    \item $\psi(\mat{A}, \vec{t}) = \psi(\mat{U}\mat{A}, \mat{U}\vec{t})$ for
      any $\mat{U} \in \mathrm{O}_{d-k}(\R)$, where $\mathrm{O}_{d-k}(\R)$ denotes the orthogonal group in dimension $(d-k)$; \label{item:range-isotropy}
    \item for $n = 1, \ldots, (d - k)$ and $m \in \N_0$,
      \begin{equation}
        \int_{\R} \psi(\mat{A}, \vec{t}) t_n^m \dd t_n
      \end{equation}
      is a degree-$m$ homogeneous polynomial with respect to $\vec{\alpha}_n$ (the $n$th row of $\mat{A}$).
  \end{enumerate}
\end{theorem}
\Cref{item:range-isotropy} is the $k$-plane analogue of the even symmetry imposed by the Radon transform. We refer to this as the \emph{isotropic symmetry} of the $k$-plane transform since it says that, under the present parameterization, the $k$-plane transform of a function is isotropic in its variables.

%% file: sections/distributional.tex
\section{Distributional Extension} \label{sec:dist-ext}

When $k=(d-1)$, there are three different ways to extend the Radon transform to distributions. These are the formulations of \begin{enumerate*}[label=(\roman*)]
    \item Gelfand, Graev, and Vilenkin~\cite{GelfandIntegralGeometry};
    \item Helgason~\cite{HelgasonIntegralGeometry}; and
    \item Ludwig~\cite{LudwigRadon}.
\end{enumerate*}
It turns out that these three approaches are equivalent~\cite{RammRadonDist}. Here, for general $k$, $1 \leq k < d$, we shall extend the $k$-plane transform to distributions by duality, taking inspiration from Ludwig~\cite{LudwigRadon}.
In \cref{sec:Fourier-slice}, we shall also make explicit that our definitions are equivalent to Fourier-based definitions understood via the Fourier slice theorem, similar to those investigated in~\cite{RammRadonDist} for the Radon transform.

In the remainder of this paper, $\ang{\dummy, \dummy}$ denotes the pairing of a space and its continuous dual defined on $\R^d$ while $\ang{\dummy, \dummy}_k$ denotes the pairing of a space and its continuous dual defined on $\Xi_k$. By invoking the property from \cref{thm:cont-inv} that $\RadonOp_k^* \KOp_{d-k} \RadonOp_k = \Id$ on $\Sch(\R^d)$ and noting that $\KOp_{d-k}$ is self-adjoint (since its frequency response is real-valued), we see that, for all $\varphi \in \Sch(\R^d)$ and $f \in \Sch'(\R^d)$,
\begin{equation}
    \ang{f, \varphi} = \ang{f, \RadonOp_k^* \KOp_{d-k} \RadonOp_k \varphi} = \ang{\RadonOp_k f, \KOp_{d-k} \RadonOp_k \varphi}_k = \ang{\RadonOp_k f, \psi}_k,
\end{equation}
where $\psi = \KOp_{d-k} \RadonOp_k\curly{\varphi} \in \KOp_{d-k} \RadonOp_k\paren{\Sch(\R^d)}$ and so $\varphi = \RadonOp_k\curly{\psi}$. Likewise, we have that
\begin{equation}
    \ang{f, \varphi} = \ang{f, \RadonOp_k^* \KOp_{d-k} \RadonOp_k \varphi} = \ang{\KOp_{d-k} \RadonOp_k f, \RadonOp_k \varphi}_k = \ang{\KOp_{d-k} \RadonOp_k f, \tilde{\psi}}_k,
\end{equation}
where $\tilde{\psi} = \RadonOp_k\curly{\varphi} \in \Sch_k$ and so $\varphi = \RadonOp_k^*\KOp_{d-k}\curly{\tilde{\psi}}$. This manipulation reveals that, in order to have a well-defined distributional extension of these operators, some care has to be taken regarding their domain and range. This is summarized in \cref{defn:distributional-k-plane}.

\begin{definition} \label[definition]{defn:distributional-k-plane} \hfill

\begin{enumerate}
    \item The \emph{distributional $k$-plane transform}
    \begin{equation}
        \RadonOp_k: \Sch'(\R^d) \to \paren*{\KOp_{d-k} \RadonOp_k\paren*{\Sch(\R^d)}}'
    \end{equation}
    is defined to be the dual map of the homeomorphism $\RadonOp_k^*: \KOp_{d-k} \RadonOp_k\paren*{\Sch(\R^d)} \to \Sch(\R^d)$.
    \label{item:distributional-k-plane}

    \item The \emph{distributional filtered $k$-plane transform}
    \begin{equation}
        \KOp_{d-k} \RadonOp_k: \Sch'(\R^d) \to \Sch_k'
    \end{equation}
    is defined to be the dual map of the homeomorphism $\RadonOp_k^* \KOp_{d-k}: \Sch_k \to \Sch(\R^d)$.
    \label{item:distributional-filtered-k-plane}

    \item The \emph{distributional backprojection}
    \begin{equation}
        \RadonOp_k^*: \Sch_k' \to \Sch'(\R^d)
    \end{equation}
    is defined to be the dual map of the homeomorphism $\RadonOp_k: \Sch(\R^d) \to \Sch_k$.

    \item The \emph{extended distributional backprojection}
    \begin{equation}
        \RadonOp_k^*: \Sch'(\Xi_k) \to \Sch'(\R^d)
    \end{equation}
    is defined to be the dual map of the continuous operator $\RadonOp_k: \Sch(\R^d) \to \Sch(\Xi_k)$, which is well-defined since $\Sch_k$ is continuously embedded in $\Sch(\Xi_k)$.
\end{enumerate}
\end{definition}
Based on these definitions, we have the following result on the invertibility of the filtered $k$-plane transform on $\Sch'(\R^d)$, which is the dual of \cref{thm:cont-inv,cor:homeo}.
\begin{theorem} \label{thm:cont-inv-dual}
    It holds that $\RadonOp_k^* \KOp_{d-k} \RadonOp_k = \Id$ on $\Sch'(\R^d)$. Moreover, the filtered $k$-plane transform $\KOp_{d-k} \RadonOp_k: \Sch'(\R^d) \to \Sch_k'$ is a homeomorphism with inverse given by the backprojection $(\KOp_{d-k} \RadonOp_k)^{-1} = \RadonOp_k^*: \Sch_k' \to \Sch'(\R^d)$.
\end{theorem}
While \cref{defn:distributional-k-plane} properly specifies the distributional $k$-plane transform, filtered $k$-plane transform, and backprojection operators, the construction is inherently abstract in the sense that the underlying spaces $\Sch_k'$ and $\paren*{\KOp_{d-k} \RadonOp_k\paren*{\Sch(\R^d)}}'$ are equivalence classes of distributions. This means that the distributional extensions of these operators lack unicity, except for the backprojection. In the case of the extended backprojection operator, we have a unique, concrete characterization of the backprojection of any distribution in $\Sch'(\Xi_k)$. Specifically, $f \in \Sch'(\R^d)$ is defined to be the backprojection of $g \in \Sch'(\Xi_k)$ (i.e., $f = \RadonOp_k^*\curly{g}$) if
\begin{equation}
    \ang{f, \varphi} = \ang{g, \RadonOp_k\curly{\varphi}}_k,
\end{equation}
for all $\varphi \in \Sch(\R^d)$. 

Conversely, the range $\Sch_k'$ of the distributional filtered $k$-plane transform can be identified with the abstract quotient space $\Sch'(\Xi_k) / \Null_{\RadonOp_k^*}$, where
\begin{equation}
    \Null_{\RadonOp_k^*} = \curly{g \in \Sch'(\Xi_k) \st \RadonOp_k^*\curly{g} = 0 \Leftrightarrow \ang{g, \tilde{\psi}}_k = 0 \text{ for all } 
    \tilde{\psi} \in \Sch_k}
    \label{eq:null-space}
\end{equation}
denotes the null space of the extended backprojection operator in \cref{defn:distributional-k-plane}. It is important to note that this null space is huge and contains many exotic functions/distributions (see~\cite{LudwigRadon} for a characterization when $k = (d-1)$).

This non-uniqueness poses an issue when one needs to deal with practical applications. To this end, in \cref{sec:k-plane-compatible} we identify Banach subspaces of $\Sch_k'$ from which we can systematically extract unique, concrete representers from the abstract equivalence classes via continuous projection operators. These \emph{$k$-plane-compatible Banach subspaces} are such that the backprojection operator is guaranteed to be invertible.

For completeness, observe that, when $k$ is even, we analogously have that the range $\paren*{\KOp_{d-k} \RadonOp_k\paren*{\Sch(\R^d)}}' = \paren*{\KOp_{d-k} \paren*{\Sch_k}}'$ of the distributional $k$-plane transform can be identified with the abstract quotient space $\Sch'(\Xi_k) / \Null_{\RadonOp_k^*\KOp_{d-k}}$, where
\begin{align*}
    &\Null_{\RadonOp_k^*\KOp_{d-k}} \\
    &= \curly{g \in \Sch'(\Xi_k) \st \RadonOp_k^*\KOp_{d-k}\curly{g} = 0 \Leftrightarrow \ang{g, \psi}_k = 0 \text{ for all } \psi \in \KOp_{d-k} \paren*{\Sch_k}}. \numberthis
\end{align*}
Here, we took advantage of the fact that, when $k$ is even, $\KOp_{d-k} = c_{d,k} (-\Delta_{d-k})^{\frac{k}{2}}$ is a local operator. In this scenario, $\KOp_{d-k} \paren*{\Sch_k}$ is continuously embedded in $\Sch(\Xi_k)$. Thus, we can define the extended distributional filtered backprojection
\begin{equation}
    \RadonOp_k^*\KOp_{d-k}: \Sch'(\Xi_k) \to \Sch'(\R^d)
\end{equation}
as the dual map of the continuous operator $\KOp_{d-k}\RadonOp_k: \Sch(\R^d) \to \Sch(\Xi_k)$. By contrast, when $k$ is odd, such a quotient-space characterization of $\paren*{\KOp_{d-k} \paren*{\Sch_k}}'$ cannot exist. Indeed, in that scenario, $\KOp_{d-k}$ is a global operator and so $\KOp_{d-k} \paren*{\Sch_k} \not\subset \Sch(\Xi_k)$.

%% file: sections/Fourier-slice.tex
\section{The Fourier Slice Theorem} \label{sec:Fourier-slice}
An important property of the $k$-plane transform is the Fourier slice theorem. It relates the Fourier transform of the $k$-plane transform of a function to its Fourier transform along \emph{slices}. While this result has been reported many times in the literature (see~\cite{MarkoeAnalyticTomo} and references therein), to the best of our knowledge, it has not yet been established for the parameterization of the $k$-plane transform given in \cref{eq:k-plane-perp,eq:k-plane-Dirac}. We first rederive the theorem for Schwartz functions with our choice of parameterization and then extend the scope of the Fourier slice theorem to tempered distributions in full generality.

\subsection{The Fourier Slice Theorem for Schwartz Functions}
In this section we begin by defining the notion of a $k$-plane ridge function. With this definition, we first prove a more general result (\cref{thm:k-plane-identity})
and then retrieve the Fourier slice theorem as a special case.
\begin{definition} \label[definition]{defn:k-plane-ridge-function}
    Given $r \in L^1_\temp(\R^{d-k})$ and $\mat{A} \in V_{d-k}(\R^d)$, the \emph{$k$-plane ridge function} $r_\mat{A} \in L^1_\temp(\R^d)$ is defined by
    \begin{equation}
        r_\mat{A}(\vec{x}) \coloneqq r(\mat{A}\vec{x}).
    \end{equation}
    The function $r$ is referred to as the \emph{profile} and the matrix $\mat{A}$ is referred to as the \emph{orientation} of the $k$-plane ridge $r_\mat{A}$, which is constant along the $k$-planes $\mat{A}\vec{x} = \vec{c}$, where $\vec{c} \in \R^{d-k}$.
\end{definition}

\begin{theorem} \label{thm:k-plane-identity}
Given $r \in L^1_\temp(\R^{d-k})$ and $\mat{A} \in V_{d-k}(\R^d)$, we have the equality
\begin{equation}
    \ang{r_\mat{A}, \varphi} = \ang{r, \RadonOp_k\curly{\varphi}(\mat{A}, \dummy)}
\end{equation}
for all $\varphi \in \Sch(\R^d)$.
\end{theorem}
\begin{proof}
    The result follows directly from the Fubini--Tonelli theorem. Indeed, we have that
    \begin{align*}
        \ang{r_\mat{A}, \varphi}
        &= \int_{\R^d} r_\mat{A}(\vec{x}) \varphi(\vec{x}) \dd\vec{x} \\
        &= \int_{\R^{d-k}} r(\mat{A}\mat{A}^\T\vec{t}) \int_{\mat{A}^\perp} \varphi(\vec{y} + \mat{A}^\T\vec{t}) \dd\vec{y} \dd\vec{t} \\
        &= \int_{\R^{d-k}} r(\vec{t}) \RadonOp_k\curly{\varphi}(\mat{A}, \vec{t}) \dd\vec{t} = \ang{r, \RadonOp_k\curly{\varphi}(\mat{A}, \dummy)} \numberthis
    \end{align*}
    under the change of variables $\vec{x} = \vec{y} + \mat{A}^\T\vec{t}$ with $\vec{y} \in \mat{A}^\perp$ and $\vec{t} \in \R^{d-k}$.
\end{proof}
Given $\mat{A} \in V_{d-k}(\R^d)$, we see from \cref{thm:k-plane-range} that $\RadonOp_k\curly{\varphi}(\mat{A}, \dummy) \in \Sch(\R^{d-k})$ for all $\varphi \in \Sch(\R^d)$. Therefore, \cref{thm:k-plane-identity} motivates the following \emph{definition} of a $k$-plane ridge whose profile does not admit a pointwise interpretation.

\begin{definition} \label[definition]{defn:k-plane-ridge-distribution}
    Given the profile $r \in \Sch'(\R^{d-k})$ and orientation $\mat{A} \in V_{d-k}(\R^d)$, the $k$-plane ridge $r_\mat{A} \in \Sch'(\R^d)$ is the tempered distribution defined by
    \begin{equation}
        \ang{r_\mat{A}, \varphi} \coloneqq \ang{r, \RadonOp_k\curly{\varphi}(\mat{A}, \dummy)}
    \end{equation}
    for all $\varphi \in \Sch(\R^d)$.
\end{definition}
\begin{remark} \label[remark]{rem:k-plane-Dirac}
    The formulation of the $k$-plane transform in \cref{eq:k-plane-Dirac} can be understood via this definition.
\end{remark}

Under the present parameterization of the $k$-plane transform, the following corollary of \cref{thm:k-plane-identity} is our reformulation of the Fourier slice theorem on $\Sch(\R^d)$.
\begin{corollary} \label[corollary]{thm:Fourier-slice-S}
    Given $\varphi \in \Sch(\R^d)$, we have, for any $\mat{A} \in V_{d-k}(\R^d)$, that
    \begin{equation}
        \FourierOp_{d-k}\curly{(\RadonOp_k\varphi)(\mat{A}, \dummy)}(\vec{\omega}) = \FourierOp_d\curly{\varphi}(\mat{A}^\T\vec{\omega}), \quad \vec{\omega} \in \R^{d-k}.
        \label{eq:Fourier-slice-S}
    \end{equation}
\end{corollary}
\begin{proof}
    Choose $r(\vec{t}) = e^{-\imag \vec{\omega}^\T\vec{t}}$ so that $r_\mat{A}(\vec{x}) = e^{-\imag \vec{\omega}^\T\mat{A}\vec{x}} = e^{-\imag (\mat{A}^\T\vec{\omega})^\T\vec{x}}$ for any $\mat{A} \in V_{d-k}(\R)$. The result then follows for any $\varphi \in \Sch(\R^d)$ by the application of \cref{thm:k-plane-identity}.    
\end{proof}

\begin{remark}
    This formulation of the Fourier slice theorem sets the stage for the investigation of the bijective isometry properties of the $k$-plane transform between $L^2$-Sobolev spaces. This has recently been established in full generality for the Radon transform ($k=(d-1)$)~\cite{RadonSobolev}. The proofs of~\cite{RadonSobolev} crucially rely on the Fourier slice theorem. This suggests that their formulation can also be extended to the $k$-plane transform.
\end{remark}

\subsection{The Fourier Slice Theorem for Tempered Distributions}
The general result for tempered distributions will be deduced by duality. To that end, we first define the operator $\FOp_{d,k}$ that acts on $\Sch(\R^d)$ by
\begin{equation}
    \FOp_{d,k}\curly{\varphi}(\mat{A}, \vec{\omega}) \coloneqq \hat{\varphi}(\mat{A}^\T\vec{\omega}) = \FourierOp_{d-k}\curly{(\RadonOp_k\varphi)(\mat{A}, \dummy)}(\vec{\omega}) = \FourierOp_{d-k}\RadonOp_k\curly{\varphi},
    \label{eq:FOp-slice}
\end{equation}
where the second equality is from \cref{thm:Fourier-slice-S}. The range of $\FOp_{d,k}$ is the space $\FourierOp_{d-k}\paren*{\Sch_k}$. The operator is invertible with its inverse being
\begin{equation}
    \FOp_{d,k}^{-1} = \RadonOp_k^* \KOp_{d-k} \FourierOp^{-1}_{d-k}.
\end{equation}
Its adjoint is given by
\begin{equation}
    \FOp_{d,k}^* = \RadonOp_k^* \FourierOp_{d-k}^* = \frac{1}{(2\pi)^{d-k}} \RadonOp_k^* \FourierOp_{d-k}^{-1}.
\end{equation}
Since $\RadonOp_k^*$ is a homeomorphism from $\KOp_{d-k} \RadonOp_k \paren*{\Sch(\R^d)} \to \Sch(\R^d)$, we have that
\begin{equation}
\FOp_{d,k}^*: \FourierOp_{d-k}\KOp_{d-k} \RadonOp_k \paren*{\Sch(\R^d)} \to \Sch(\R^d)
\end{equation}
is a homeomorphism. This justifies the following extension of $\FOp_{d,k}$ by duality:
\begin{equation}
    \FOp_{d,k}: \Sch'(\R^d) \to \paren*{\FourierOp_{d-k}\KOp_{d-k} \RadonOp_k \paren*{\Sch(\R^d)}}'
    \label{eq:FOp-dist}
\end{equation}
which is itself a homeomorphism. Next, from \cref{eq:FOp-slice}, we have the identity
\begin{equation}
    \RadonOp_k = \FourierOp_{d-k}^{-1} \FOp_{d,k} \quad\text{on } \Sch(\R^d).
    \label{eq:k-plane-Fourier}
\end{equation}
The combination of \cref{eq:k-plane-Fourier} with \cref{thm:cont-inv} yields the identity
\begin{equation}
    \KOp_{d-k} \RadonOp_k = \FourierOp_{d-k}^{-1} \hat{K}_{d-k} \FOp_{d,k} \quad\text{on } \Sch(\R^d),
\end{equation}
where $\hat{K}_{d-k}$ denotes multiplication in the $\vec{\omega}$ variable by the frequency response $\hat{K}_{d-k}(\vec{\omega}) = c_{d,k} \norm{\vec{\omega}}_2^k$ of the filtering operator $\KOp_{d-k}$. From \cref{thm:cont-inv}, $\paren{\RadonOp_k^*}^{-1} = \KOp_{d-k} \RadonOp_k$ on $\Sch(\R^d)$. Therefore, this also establishes the identity
\begin{equation}
    \RadonOp^* = \FOp_{d,k}^{-1} \hat{K}_{d-k}^{-1} \FourierOp_{d-k} \quad \text{on } \KOp_{d-k} \RadonOp_k \paren*{\Sch(\R^d)},
    \label{eq:backprojection-Fourier}
\end{equation}
where $\hat{K}_{d-k}^{-1}$ denotes multiplication in the $\vec{\omega}$ variable by the frequency response $\hat{K}_{d-k}^{-1}(\vec{\omega}) = c_{d,k}^{-1} \norm{\vec{\omega}}_2^{-k}$ of the inverse filtering operator $\KOp_{d-k}^{-1}$. 

\Cref{eq:backprojection-Fourier} is a new formula for the dual transform. In particular, if we fix $f \in \Sch'(\R^d)$, for any $\psi \in \KOp_{d-k} \RadonOp_k \paren*{\Sch(\R^d)}$ we have
\begin{equation}
    \ang{\RadonOp_k f, \psi}_k = \ang{f, \FOp_{d,k}^{-1} \hat{K}_{d-k}^{-1} \FourierOp_{d-k} \psi}
\end{equation}
By noticing that $\FOp_{d,k}^{-1}\hat{K}_{d-k}^{-1} = \RadonOp_k^* \KOp_{d-k} \FourierOp^{-1}_{d-k} \hat{K}_{d-k}^{-1} = \RadonOp_k^* \FourierOp^{-1}_{d-k} = (2\pi)^{d-k} \FOp_{d,k}^*$, we have that
\begin{equation}
    \ang{\RadonOp_k f, \psi}_k = (2\pi)^{d-k} \ang{\FOp_{d,k} f, \FourierOp_{d-k} \psi}_k,
    \label{eq:Fourier-slice-precursor}
\end{equation}
for any $\psi \in \KOp_{d-k} \RadonOp_k \paren*{\Sch(\R^d)}$, where the right-hand side is well-defined from \cref{eq:FOp-dist}. An immediate consequence to \cref{eq:Fourier-slice-precursor} is the Fourier slice theorem for any $f \in \Sch'(\R^d)$.
\begin{theorem} \label{thm:Fourier-slice}
    Let $f \in \Sch'(\R^d)$. Then,
    \begin{equation}
        \FourierOp_{d-k} \RadonOp_k f = \FOp_{d,k} f,
        \label{eq:Fourier-slice-equality}
    \end{equation}
    where the equality holds in $\paren*{\FourierOp_{d-k}\KOp_{d-k} \RadonOp_k \paren*{\Sch(\R^d)}}'$.
\end{theorem}
When $f \in \Sch(\R^d)$, the equality in \cref{eq:Fourier-slice-equality} is precisely the equality in \cref{eq:Fourier-slice-S} and so \cref{thm:Fourier-slice} recovers \cref{thm:Fourier-slice-S}.
For general tempered distributions, the identity in \cref{thm:Fourier-slice} has to be interpreted with care since
$\paren*{\KOp_{d-k} \RadonOp_k \paren*{\Sch(\R^d)}}'$ is actually a space of equivalence classes of distributions.
This is to say that, given $f \in \Sch'(\R^d)$, $\RadonOp_k\curly{f}$ is, in general, non-unique. This is illustrated concretely in~\cite[Example~10.3.2]{RammRadonBook}, which provides an explicit calculation via the Fourier slice theorem to show that the Radon transform ($k = (d-1)$) of the function $f \equiv 1 \in \Sch'(\R^d)$ is not uniquely determined.

The key takeaway message from the Fourier slice theorem (\cref{thm:Fourier-slice-S,thm:Fourier-slice}) is that the restriction of $\FourierOp_d\curly{f}$ to the slice $\curly{\mat{A}^\T\vec{\omega} \in \R^d \st \vec{\omega} \in \R^{d-k}}$ is the $(d-k)$-variate Fourier transform of $\RadonOp_k\curly{f}(\mat{A}, \dummy)$. Moreover, a further consequence of this theorem is a general form of the \emph{intertwining properties} of the $k$-plane transform. Before we state our result, recall that the necessary and sufficient condition for an LSI operator $\LOp_\rad$ to continuously map $\Sch(\R) \to \Sch(\R)$ (or, by duality, continuously map $\Sch'(\R) \to \Sch'(\R)$) is that $\hat{L}_\rad \in \mathcal{O}_M(\R)$, where $\mathcal{O}_M(\R)$ is the space of smooth functions of slow growth on $\R$~\cite[p.~243]{Schwartz}. 

\begin{corollary} \label[corollary]{thm:intertwine}
    Let $\LOp_{d-k}: \Sch'(\R^{d-k}) \to \Sch'(\R^{d-k})$ and $\LOp_d: \Sch'(\R^d) \to \Sch'(\R^d)$ be isotropic LSI operators that share the same radial frequency profile. Their frequency responses satisfy
    \begin{equation}
        \hat{L}_{d-k}(\vec{\omega}) = \hat{L}_\rad(\norm{\vec{\omega}}_2), \quad \vec{\omega} \in \R^{d-k},
    \end{equation}
    and
    \begin{equation}
        \hat{L}_d(\vec{\xi}) = \hat{L}_\rad(\norm{\vec{\xi}}_2), \quad \vec{\xi} \in \R^d,
    \end{equation}
    where $\hat{L}_\rad \in \mathcal{O}_M(\R)$ is their shared radial frequency profile. Then, we have the intertwining properties
    \begin{equation}
        \LOp_{d-k} \RadonOp_k = \RadonOp_k \LOp_d \quad \text{on } \Sch'(\R^d)
    \label{eq:intertwining-thm}
    \end{equation}
    and
    \begin{equation}
      \RadonOp_k^* \LOp_{d-k}^* = \LOp_d^* \RadonOp_k^* \quad \text{on } \Sch_k'.
    \end{equation}
\end{corollary}
\begin{proof}
    With a slight abuse of notation, by \cref{thm:Fourier-slice}, one has that
\begin{align*}
    \FourierOp_{d-k} \LOp_{d-k} \RadonOp_k
    &= \hat{L}_{d-k}(\vec{\omega}) \FourierOp_{d-k} \RadonOp_k \\
    &= \hat{L}_\rad(\norm{\vec{\omega}}_2) \FOp_{d,k} \\
    &= \hat{L}_\rad(\norm{\mat{A}^\T\vec{\omega}}_2) \FOp_{d,k} \\
    &= \hat{L}_d(\mat{A}^\T\vec{\omega}) \FOp_{d,k} \\
    &= \FOp_{d,k} \LOp_d \\
    &= \FourierOp_{d-k} \RadonOp_k \LOp_d. \numberthis
\end{align*}
Furthermore, since \cref{eq:intertwining-thm} holds on $\Sch(\R^d)$, by duality, we immediately have an intertwining property for the dual transform on $\Sch_k'$ as
\begin{equation}
  \RadonOp_k^* \LOp_{d-k}^* = \LOp_d^* \RadonOp_k^* \quad \text{on } \Sch_k'.
\end{equation}
\end{proof}

%% file: sections/compatible.tex
\section{Systematic Construction of \texorpdfstring{$k$}{k}-Plane-Compatible Spaces} \label{sec:k-plane-compatible}
The distributional extension of the $k$-plane transform and related operators in \cref{sec:dist-ext} provides a definition of the filtered $k$-plane transform of any $f \in \Sch'(\R^d)$. It also provides a way to invert it via the backprojection operator. So far, this is not very practical since the range of the distributional filtered $k$-plane transform $\Sch_k'$ is a space of equivalence classes of distributions which has been identified as the quotient $\Sch'(\Xi_k) / \Null_{\RadonOp_k^*}$. To make the framework more applicable, we now provide a systematic construction (via a completion process) of Banach subspaces of $\Sch_k'$ for which the backprojection is invertible and whose elements are identified as ordinary functions or measures. We shall refer to these as $k$-plane-compatible (Banach) spaces. We then present some examples where, in particular, we can identify unique, concrete representers from the equivalence classes in $\Sch_k' \cong \Sch'(\Xi_k) / \Null_{\RadonOp_k^*}$ via continuous projection operators. 

Let $\X$ be a Banach space on $\Xi_k$ such that
\begin{equation}
  \Sch(\Xi_k) \dembed \X \dembed \Sch'(\Xi_k),
\end{equation}
the corresponding dual pair of $k$-plane-compatible Banach spaces being
\begin{equation}
  \paren{\X_k, \X_k'},
\end{equation}
with $\X_k \coloneqq \cl{(\Sch_k, \norm*{\dummy}_\X)}$, and for which we have the following powerful result.

\begin{theorem} \label{thm:k-plane-systematic}
  Let $(\X_k, \X_k')$ be as above. Then, the following holds.
  \begin{enumerate}
    \item The backprojection $\RadonOp_k^*: \X_k' \to \Sch'(\R^d)$ is injective
      and $\KOp_{d-k}\RadonOp_k \RadonOp_k^* = \Id$ on $\X_k'$.
      \label{item:range}

    \item Let $\Y = \cl{(\Sch(\R^d), \norm{\dummy}_\Y)}$ where, for $\varphi \in
      \Sch(\R^d)$, $\norm{\varphi}_\Y \coloneqq \norm{\RadonOp_k \varphi}_\X$.
      Then, there exists a unique extension $\RadonOp_k^*\KOp_{d-k}: \X_k \to \Y$ that 
      is an isometric isomorphism such that $\RadonOp_k \RadonOp_k^* \KOp_{d-k} =
      \Id$ on $\X_k$. \label{item:X_k-Y-iso}

    \item The range spaces $\Y = \RadonOp_k^*\KOp_{d-k}\paren{\X_k}$ and $\V =
      \RadonOp_k^*\paren{\X_k'}$ form a dual pair of Banach spaces (i.e., we can
      identify $\Y' = \V$), such that $(\Y, \V)$ is isometrically isomorphic to
      $(\X_k, \X_k')$.
  \end{enumerate}
  Moreover, if there exists a complementary Banach space $\X_k^\comp$ such that $\X = \X_k \oplus \X_k^\comp$, then the following holds.
  \begin{enumerate}[resume]
      \item The dual space is decomposable as $\X' = \X_k' \oplus \paren{\X_k^\comp}'$. \label{item:direct-sum-dual-Xk}
      \item The dual complement space $\paren{\X_k^\comp}'$ is the null space of $\RadonOp_k^*: \X' \to \Y'$. \label{item:dual-complement-null-space}
      \item The complement space $\X_k^\comp$ is the null space of $\RadonOp_k^* \KOp_{d-k}: \X \to \Y$. \label{item:complement-null-space}
      \item The operators $\P_k = \RadonOp_k \RadonOp_k^* \KOp_{d-k}: \X \to \X_k$ and $\P_k^* = \KOp_{d-k}\RadonOp_k \RadonOp_k^*: \X' \to \X_k'$ form an adjoint pair of continuous projectors with $\P_k\paren{\X} = \X_k$ and $\P_k^*\paren{\X'} = \X_k'$.
  \end{enumerate}
\end{theorem}
\begin{proof} \hfill

  \begin{enumerate}
    \item Since $\Sch_k \dembed \X_k$, by duality we have that $\X_k' \cembed \Sch_k'$. Let
      $\V = \RadonOp_k^*\paren{\X_k'} \subset \Sch'(\R^d)$ denote the
      range of the restricted operator $\eval{\RadonOp_k^*}_{\X_k'}$. Clearly
      this restriction is a homeomorphism from $\X_k' \to \V$, where 
      $\V$ is a Banach space equipped with the norm
      \begin{equation}
        \norm{f}_\V \coloneqq
        \norm{\RadonOp_k^{*-1}\curly{f}}_{\X_k'}, \quad f \in \V.
      \end{equation}
      In other words, $\RadonOp_k^*$ is an isometric isomorphism from $\X_k'$ to $\V$.
      In particular, this shows that $\RadonOp_k^*: \X_k' \to \Sch'(\R^d)$ is an injection.
      Moreover, by the inversion formula in \cref{thm:cont-inv-dual}, the inverse of
      $\eval{\RadonOp_k^*}_{\X_k'}$ on $\V$ is given by the restriction of
      $\KOp_{d-k}\RadonOp_k$ to $\V$, denoted $\eval{\KOp_{d-k}\RadonOp_k}_\V$.
      Therefore, $\KOp_{d-k}\RadonOp_k\RadonOp_k^* = \Id$ on $\X_k'$.

    \item From the proof of \cref{item:range}, we have that $\KOp_{d-k} \RadonOp_k: \V \to \X_k'$ is an isometric isomorphism and, hence, it is continuous by design. By duality combined with the isometric embedding of a Banach space in its bidual, we have that $\RadonOp_k^* \KOp_{d-k}: \X_k \to \V'$ is continuous as well. We now prove that it is injective by showing that it has a trivial null space. For any $\psi \in \Sch_k$, the statement $\RadonOp_k^*\KOp_{d-k}\curly{\psi} = 0$ is equivalent to
    \begin{equation}
        0 = \ang{v, \RadonOp_k^*\KOp_{d-k}\curly{\psi}} = \ang{\KOp_{d-k} \RadonOp_k \curly{v}, \psi}_k, \quad  \text{for all } v \in \V.
        \label{eq:RK-null-trivial}
    \end{equation}
    Since $\KOp_{d-k}\RadonOp_k\paren{\V} = \X_k'$, \cref{eq:RK-null-trivial} is equivalent to $\ang{g, \psi}_k = 0$ for all $g \in \X_k'$, which implies $\psi = 0$. The result then follows from the dense embeddings $\Sch_k \dembed
      \X_k$ and $\Sch(\R^d) \dembed \Y$ and the continuity of the inversion
      formula for the $k$-plane transform. Indeed, these dense embeddings imply that the definition of $\Y$ is equivalent to
      \begin{equation}
        \Y = \curly{\RadonOp_k^{-1}\curly{g} \st g \in \X_k},
      \end{equation}
      where $\RadonOp_k^{-1} = \RadonOp_k^*\KOp_{d-k}$ is continuous on $\Sch_k$ and
      can can be continuously extended, by density, to act on $\X_k$. Therefore,
      $\Y = \RadonOp_k^*\KOp_{d-k}\paren{\X_k}$. Since $\RadonOp_k^*\KOp_{d-k}$ is continuous and injective on $\X_k$, by the bounded inverse theorem~\cite{RudinFA}, it has a continuous inverse on its range that maps $\Y \to \X_k$. Putting everything together, from the inversion formula for
      the $k$-plane transform and the definition of the norm of $\Y$, we have
      that $\RadonOp_k^*\KOp_{d-k}: \X_k \to \Y$ is an isometric isomorphism with
      inverse given by $(\RadonOp_k^*\KOp_{d-k})^{-1} = \RadonOp_k: \Y \to \X_k$ and,
      so, $\RadonOp_k \RadonOp_k^* \KOp_{d-k} = \Id$ on $\X_k$.

    \item From \cref{item:X_k-Y-iso}, we have that
      \begin{align*}
        \RadonOp_k^*\KOp_{d-k}: \X_k \to \Y \\
        \RadonOp_k: \Y \to \X_k \numberthis
      \end{align*}
      are bijective isometries. So, by duality, we have that
      \begin{align*}
        \KOp_{d-k} \RadonOp_k: \Y' \to \X_k' \\
        \RadonOp_k^*: \X_k' \to \Y' \numberthis
      \end{align*}
      are bijective isometries. From the
      proof of \cref{item:range}, we know that $\X_k'$ is isometrically
      isomorphic to $\V = \RadonOp_k^*\paren{\X_k'}$. This means that $\V$ is
      isometrically isomorphic to $\Y'$. Therefore, we can identify $\Y' = \V$.

      \item This follows from the generic property that the dual of a direct sum is isometrically isomorphic to the direct sum of the duals.
      
      \item \Cref{item:direct-sum-dual-Xk} implies that $\paren{\X_k^\comp}'$ is the annihilator of
      $\X_k$ in $\X'$ and so
      \begin{align*}
        \paren{\X_k^\comp}' 
        &= \curly{g \in \X' \st \ang{g, \psi}_k = 0 \text{ for all } \psi \in \X_k} \\
        &=\curly{g \in \X' \st \ang{g, \psi}_k = 0 \text{ for all } \psi \in \Sch_k}, \numberthis
      \end{align*}
      where the second line follows from the denseness of $\Sch_k$ in $\X_k$. From the definition of the null space in \cref{eq:null-space}, we see that $\paren{\X_k^\comp}' \subset \Null_{\RadonOp_k^*}$ is the null space of $\RadonOp_k^*: \X' \to \Y'$.

      \item This result follows from \cref{item:dual-complement-null-space}, by duality.

      \item We first recall that existence (unicity) and continuity of the projectors are guaranteed for any pair of complemented Banach spaces. Next, \Cref{item:range,item:X_k-Y-iso} imply that $\P_k^*\paren{\X_k'} = \X_k'$ and $\P_k\paren{\X_k} = \X_k$. Similarly, \Cref{item:dual-complement-null-space,item:complement-null-space} imply that $\P_k^*\paren{\paren{\X_k^\comp}'} = \curly{0}$ and $\P_k\paren{\X_k^\comp} = \curly{0}$. Therefore, $(\P_k, \P_k^*)$ form an adjoint pair of continuous projectors with $\P_k\paren{\X} = \X_k$ and $\P_k^*\paren{\X'} = \X_k'$.
  \end{enumerate}
\end{proof}

\begin{figure}[t!]
    \centering
    \includegraphics[width=0.7\textwidth]{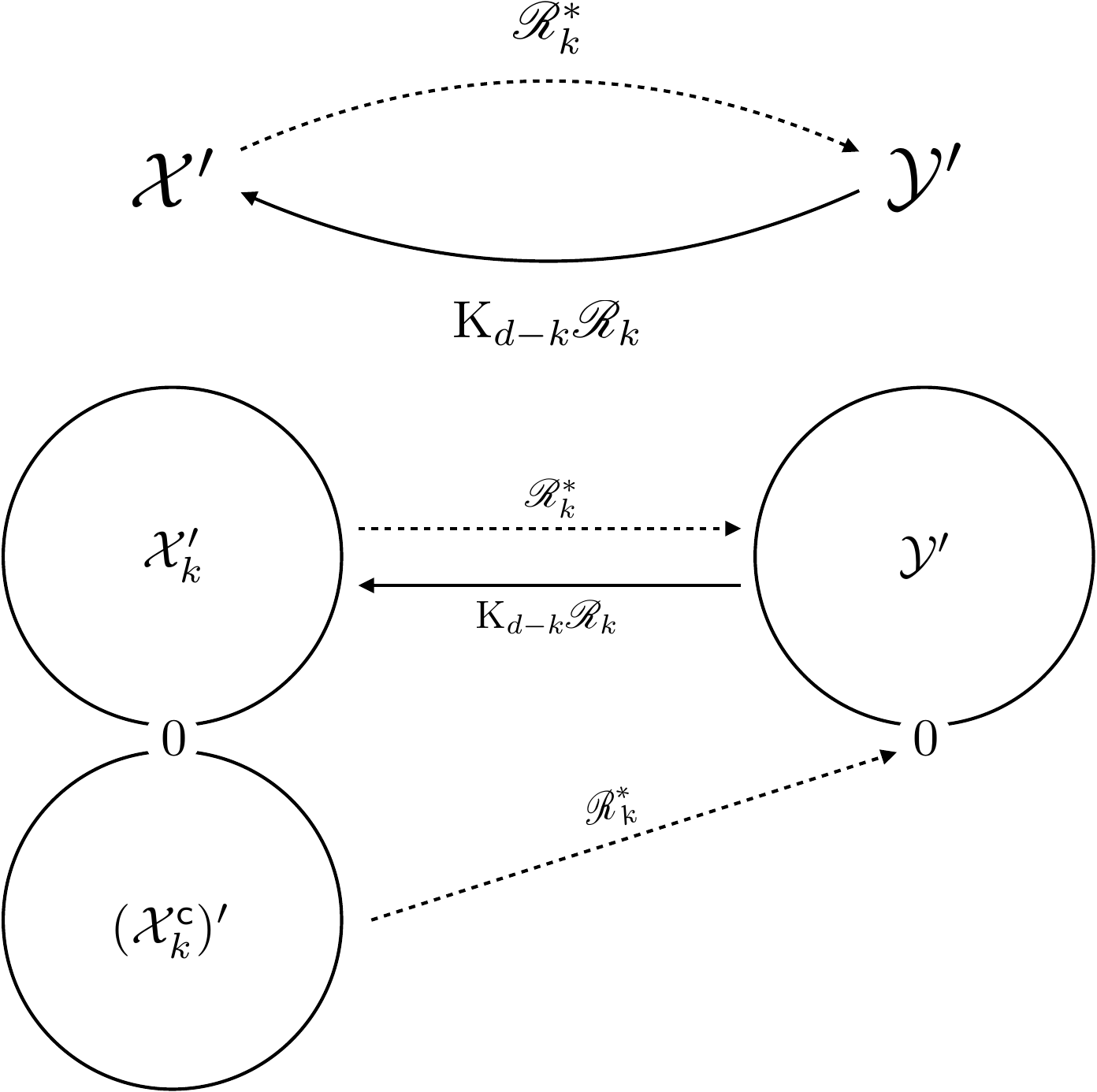}
    \caption{Relevant Banach spaces when $\X = \X_k \oplus \X_k^\comp$.}
    \label{fig:spaces}
\end{figure}
When there exists a complementary Banach space $\X_k^\comp$ such that $\X = \X_k \oplus \X_k^\comp$, we depict the relevant Banach spaces from \cref{thm:k-plane-systematic} and their relations in \cref{fig:spaces}. Our next theorem and corollary characterizes some useful examples of $k$-plane-compatible Banach spaces whose elements are identifiable as concrete ``isotropic'' functions or measures (from  the abstract equivalence classes in $\Sch_k' \cong \Sch'(\Xi_k) / \Null_{\RadonOp_k}$).

\begin{theorem} \label{thm:k-plane-compatible}
    Let $\alpha_p$ denote the natural tensor norm for tensor products of $L^p$ spaces and let $\otimes_{\alpha_p}$ denote the completion of the algebraic tensor product with respect to the tensor norm $\alpha_p$. If $\X$ is any of the following Banach spaces:
    \begin{itemize}
        \item $L^p(\Xi_k) = L^p(V_{d-k}(\R^d) \otimes_{\alpha_p} L^p(\R^{d-k})$ for $1 < p < \infty$;
        
        \item $C_0(\Xi_k) = C(V_{d-k}(\R^d)) \otimes_{\alpha_\infty} C_0(\R^{d-k})$,
    \end{itemize}
    then $\X = \X_k \oplus \X_k^\comp$, where $\X_k = \X_\iso$, where
    \begin{equation}
        \X_\iso \coloneqq \curly{g \in \X \st g(\mat{A}, \vec{t}) = g(\mat{U}\mat{A}, \mat{U}\vec{t}), \text{ for any } (\mat{A}, \vec{t}) \in \Xi_k \text{ and } \mat{U} \in \rmO_{d-k}(\R)}
    \end{equation}
    is the subspace of isotropic functions in $\X$.
\end{theorem}
\begin{proof}
    The space $\X_\iso$ is a closed subspace of $\X$ by construction. To show that it is complemented in $\X$, we consider the self-adjoint operator $\P_\iso: \X \to \X$ (see \cref{app:iso-projector}) which extracts the isotropic part of a function as
    \begin{equation}
        \P_\iso\curly{g}(\mat{A}, \vec{t}) = \avint_{\mathrm{O}_{d-k}(\R)} g(\mat{U}\mat{A}, \mat{U}\vec{t}) \dd\sigma(\mat{U}),
        \label{eq:P-iso-proof}
    \end{equation}
    where $\avint$ denotes the average integral and $\sigma$ is the Haar measure on $\rmO_{d-k}(\R)$. Clearly, $\P_\iso(\X) \subset \X_\iso$ and, given $g \in \X_\iso$, $\P_\iso\P_\iso g = g$. This establishes that $\P_\iso$ is a projector and that $\P_\iso(\X) = \X_\iso$. This projector is bounded since
    \begin{align*}
        \norm{\P_\iso g}_\X
        &= \norm*{(\mat{A}, \vec{t}) \mapsto \avint_{\mathrm{O}_{d-k}(\R)} g(\mat{U}\mat{A}, \mat{U}\vec{t}) \dd\sigma(\mat{U})}_\X \\
        &\leq \avint_{\mathrm{O}_{d-k}(\R)} \norm*{(\mat{A}, \vec{t}) \mapsto g(\mat{U}\mat{A}, \mat{U}\vec{t})}_\X  \dd\sigma(\mat{U}) \\
        &= \avint_{\mathrm{O}_{d-k}(\R)} \norm{g}_\X  \dd\sigma(\mat{U}) \\
        &= \norm{g}_\X, \numberthis
    \end{align*}
    where the second line follows from Jensen's inequality and the third line follows from the rotational invariance of the $L^p$-norms. The boundedness of the projector $\P_\iso$ guarantees the existence of a complementary projector $(\Id - \P_\iso): \X \to \X$ whose (closed) range $(\Id - \P_\iso)(\X) \eqqcolon \X_\iso^\comp$ is the complementary subspace of $\X_\iso$ with the properties that $\X = \X_\iso \oplus \X_\iso^\comp$ and so $\X_\iso \cap \X_\iso^\comp = \curly{0}$.

    The final part of the proof is to show that
    \begin{equation}
        \X_k = \cl{(\Sch_k, \norm{\dummy}_\X)} = \X_\iso.
        \label{eq:Xk-Xiso}
    \end{equation}
    To that end, we consider the Lizorkin space
    \begin{equation}
        \Sch_\Liz(\R^n) = \curly*{\varphi \in \Sch(\R^n) \st \int_{\R^n} \vec{x}^\vec{m} \varphi(\vec{x}) \dd\vec{x} = 0, \text{ for every multi-index } \vec{m}},
    \end{equation}
    which is the closed subspace of Schwartz functions with all moments vanishing. We can define the Lizorkin space on $\Xi_k$ as $\Sch_\Liz(\Xi_k) \coloneqq C^\infty(V_{d-k}(\R^d)) \,\hat{\otimes}\, \Sch_\Liz(\R^{d-k})$. Then, let
    \begin{equation}
        \Sch_{\Liz, \iso}(\Xi_k) = \curly*{\psi \in \Sch_\Liz(\Xi_k) \st
        \begin{aligned}
        &\psi(\mat{A}, \vec{t}) = \psi(\mat{U}\mat{A}, \mat{U}\vec{t}) \\
        &\text{for any } (\mat{A}, \vec{t}) \in \Xi_k \text{ and } \mat{U} \in \rmO_{d-k}(\R)
        \end{aligned}}
    \end{equation}
    denote its isotropic subspace. From \cref{thm:k-plane-range}, we see that $\Sch_{\Liz, \iso}(\Xi_k) \subset \Sch_k$. 

    Since
    \begin{equation}
        \cl{(C^\infty(V_{d-k}(\R^d)), \norm{\dummy}_{L^p})} = \begin{cases}
            L^p(V_{d-k}(\R^d)), & 1 < p < \infty \\
            C(V_{d-k}(\R^d)), & p = \infty
        \end{cases}
        \label{eq:cl1}
    \end{equation}
    and
    \begin{equation}
        \cl{(\Sch_\Liz(\R^{d-k}), \norm{\dummy}_{L^p})} = \begin{cases}
            L^p(\R^{d-k}), & 1 < p < \infty \text{ (from~\cite{SamkoDense})} \\
            C_0(\R^{d-k}), & p = \infty \text{ (from~\cite{NeumayerLizorkin})},
        \end{cases}
        \label{eq:cl2}
    \end{equation}
    and since the isotropic subspace is complete, then \cref{eq:cl1,eq:cl2} imply that $\X_\iso = \cl{(\Sch_{\Liz, \iso}(\Xi_k), \norm{\dummy}_\X)}$. Since $\Sch_k \subset \X_\iso$, this proves \cref{eq:Xk-Xiso}, since the completion is unique.
\end{proof}
For these $k$-plane-compatible Banach spaces, the proof of \cref{thm:k-plane-compatible} actually implies that $\P_k = \P_k^* = \P_\iso$. In other words, it is the projection operator that extracts the isotropic part of a function (or distribution).
From \cref{thm:k-plane-systematic} we also know that the backprojection operator is invertible on $\X_\iso'$, where $\X$ is any of the Banach spaces in \cref{thm:k-plane-compatible}. This is summarized in \cref{cor:backprojection-invertible}.

\begin{corollary} \label[corollary]{cor:backprojection-invertible}
    The backprojection operator $\RadonOp_k^*$ is invertible on
    \begin{enumerate}
        \item $\M_\iso(\Xi_k)$; \label{item:M-invert}
        \item $L^p_\iso(\Xi_k)$ for $1 < p < \infty$. \label{item:Lp-invert}
    \end{enumerate}
\end{corollary}
\begin{proof}
    This follows by duality directly from \cref{thm:k-plane-compatible}.
\end{proof}

\begin{remark} \label[remark]{rem:dual-Radon}
    \Cref{thm:k-plane-systematic} recovers the result of~\cite[Theorem~8]{UnserRidges} in the special case of $k = (d-1)$. To the best of our knowledge, not much was previously known about the invertibility of $\RadonOp_k^*$ beyond the fact that it is invertible on $\Sch_k'$. However, in the special case $k = (d-1)$, a multitude of classical works have addressed the invertibility of the dual Radon transform $\RadonOp_{d-1}^*$. In particular, it is known that it is invertible on
    \begin{enumerate}
        \item the subspace $\Sch_\even(\cyl)$ of even functions of the Schwartz space of smooth and rapidly decreasing functions on $\cyl$~\cite{SolmonDualRadon};
        
        \item the space of essentially bounded, even, and compactly supported functions on $\cyl$, denoted as ~$L^\infty_{\even, \cpct}(\cyl)$~\cite{RammBackprojection};
    
        \item the subspace $\Sch_{\Liz, \even}'(\cyl)$ of even Lizorkin distributions on $\cyl$~\cite{KostadinovaRidgelets}.
    \end{enumerate}
    Therefore, \cref{cor:backprojection-invertible} is among the first to provide results on the invertibility of $\RadonOp_k^*$, some of which are even new for the dual Radon transform, notably, the invertibility of $\RadonOp_{d-1}^*$ on $L^p_\even(\cyl)$ for $1 < p < \infty$. In particular, since both $\Sch_\iso(\Xi_k)$ and $L^\infty_{\iso, \cpct}(\Xi_k)$ continuously embed into $L^p_\iso(\Xi_k)$ for any $1 < p < \infty$, we have that \cref{cor:backprojection-invertible} provides a direct generalization of the results of~\cite{RammBackprojection,SolmonDualRadon} for the $k$-plane transform.
\end{remark}

%% file: sections/specific-transforms.tex
\section{Specific \texorpdfstring{$k$}{k}-Plane Transforms} \label{sec:specific-transforms}
In this section, we provide explicit formulas for the $k$-plane transform of isotropic functions and $k$-plane ridge distributions.

\subsection{Isotropic Functions}
An isotropic function $\rho_\iso: \R^d \to \R$ is characterized by its radial profile $\rho: \R \to \R$ via the equality $\rho_\iso(\vec{x}) = \rho(t)$, where $t = \norm{\vec{x}}_2$. The frequency-domain counterpart of this characterization is $\hat{\rho}_\iso(\vec{\xi}) = \hat{\rho}_\rad(\omega)$, where $\omega = \norm{\vec{\xi}}_2$. The radial frequency profile $\hat{\rho}_\rad$ can be computed via the Hankel transform~\cite[Theorem~5.26]{WendlandBook}
\begin{equation}
    \hat{\rho}_\rad(\omega) = \frac{(2\pi)^{d/2}}{\abs{\omega}^{d/2 - 1}} \int_0^\infty J_{d/2-1}(\omega t) t^{d/2 - 1} \rho(t) t \dd t,
\end{equation}
where $J_\nu$ is is the Bessel function of the first kind of order $\nu$. 

\begin{theorem} \label{thm:k-plane-isotropic}
    Let $\rho_\iso \in \Sch(\R^d)$ be an isotropic function with radial frequency profile $\hat{\rho}_\rad \in \Sch(\R)$. Then, for any $\vec{x}_0 \in \R^d$,
    \begin{equation}
        \RadonOp_k\curly{\rho_\iso(\dummy - \vec{x}_0)}(\mat{A}, \vec{t}) = \tilde{\rho}_\rad(\vec{t} - \mat{A}\vec{x}_0),
    \end{equation}
    where $\tilde{\rho}_\rad = \FourierOp_{d-k}^{-1}\curly{\hat{\rho}_\rad(\norm{\vec{\dummy}}_2)}$ and $\FourierOp^{-1}_{d-k}$ denotes the $(d-k)$-variate inverse Fourier transform.
\end{theorem}
\begin{proof}
    By the Fourier slice theorem,
    \begin{equation}
        \FourierOp_{d-k}\curly{\RadonOp_k\curly{\rho_\iso}(\mat{A}, \dummy)}(\vec{\omega}) = \hat{\rho}_\iso(\mat{A}^\T\vec{\omega}) = \hat{\rho}_\rad(\norm{\mat{A}^\T\vec{\omega}}_2) = \hat{\rho}_\rad(\norm{\vec{\omega}}_2).
        \label{eq:Fourier-iso}
    \end{equation}
    The $(d-k)$-variate inverse Fourier transform of \cref{eq:Fourier-iso} yields that
    \begin{equation}
        \RadonOp_k\curly{\rho_\iso}(\mat{A}, \vec{t}) = \tilde{\rho}_\rad(\vec{t}).
    \end{equation}
    The result then follows by the shifting property of the $k$-plane transform. Indeed, for any $\varphi \in \Sch(\R^d)$ we have that
    \begin{align*}
        \RadonOp_k\curly{\varphi(\dummy - \vec{x}_0)}(\mat{A}, \vec{t})
        &= \int_{\R^d} \varphi(\vec{x} - \vec{x}_0) \delta(\mat{A}\vec{x} - \vec{t}) \dd\vec{x} \\
        &= \int_{\R^d} \varphi(\vec{y}) \delta(\mat{A}\vec{y} - (\vec{t} - \mat{A}\vec{x}_0)) \dd\vec{y} \\
        &= \RadonOp_k\curly{\varphi}(\mat{A}, \vec{t} - \mat{A}\vec{x}_0). \numberthis
    \end{align*}
\end{proof}
\begin{remark}
    \Cref{thm:k-plane-isotropic} can be extended to general isotropic tempered distributions so long as the equality is understood in the distribution space $\paren*{\KOp_{d-k} \RadonOp_k \paren*{\Sch(\R^d)}}'$.
\end{remark}

\begin{example} \label[example]{example:k-plane-Gaussian}
    Consider the isotropic Gaussian density with zero mean and unit variance
    \begin{equation}
        g_\iso(\vec{x}) = \frac{1}{(2\pi)^{d/2}} \exp\paren*{-\frac{\norm{\vec{x}}_2^2}{2}}.
    \end{equation}
    Its Fourier transform is given by $\hat{g}_\iso(\vec{\xi}) = \exp\paren*{-\norm{\vec{\xi}}_2^2 / 2}$, which corresponds to the radial profile
    $\hat{g}_\rad(\omega) = \exp\paren*{-\abs{\omega}^2 / 2}$. The determination of the $(d-k)$-variate inverse Fourier transform then gives
    \begin{equation}
        \FourierOp_{d-k}^{-1}\curly{\hat{g}_\rad(\norm{\vec{\dummy}}_2)}(\vec{t}) = \frac{1}{(2\pi)^{\frac{d-k}{2}}} \exp\paren*{-\frac{\norm{\vec{t}}_2^2}{2}}.
    \end{equation}
    Finally, by applying \cref{thm:k-plane-isotropic}, we find for any $\vec{x}_0 \in \R^d$, that
    \begin{equation}
        \RadonOp_k\curly{g_\iso(\dummy - \vec{x}_0)}(\mat{A}, \vec{t}) = \frac{1}{(2\pi)^{\frac{d-k}{2}}} \exp\paren*{-\frac{\norm{\vec{t} - \mat{A}\vec{x}_0}_2^2}{2}}.
    \end{equation}
\end{example}

\subsection{\texorpdfstring{$k$}{k}-Plane Ridges}
From \cref{defn:k-plane-ridge-distribution}, we have the definition of of a $k$-plane ridge with any profile $r \in \Sch'(\R^{d-k})$. In this section, we write $r(\mat{A}_0(\dummy) - \vec{t}_0) \in \Sch'(\R^d)$ to denote the $k$-plane ridge with profile $r(\dummy - \vec{t}_0) \in \Sch'(\R^{d-k})$, where $(\mat{A}_0, \vec{t}_0) \in \Xi_k$.

\begin{theorem} \label{thm:k-plane-ridges}
    Let $(\mat{A}_0, \vec{t}_0) \in \Xi_k$ and $r \in \Sch'(\R^{d-k})$. Then,
    \begin{align}
        \KOp_{d-k}\RadonOp_k\curly*{r(\mat{A}_0(\dummy) - \vec{t}_0)} &= \sq*{\delta(\dummy - \mat{A}_0) \, r(\dummy - \vec{t}_0)} \in \Sch_k' \label{eq:filtered-k-plane-ridge-abstract} \\ 
        \RadonOp_k\curly*{r(\mat{A}_0(\dummy) - \vec{t}_0)} &= \sq*{\delta(\dummy - \mat{A}_0) \, (q * r)(\dummy - \vec{t}_0)} \in \paren*{\KOp_{d-k} \RadonOp_k\paren*{\Sch(\R^d)}}', \label{eq:k-plane-ridge-abstract}
    \end{align}
    where $q(\vec{t}) = \FourierOp_{d-k}^{-1}\curly{\hat{K}_k^{-1}}(\vec{t})$ is the impulse response of the inverse filtering operator $\KOp_{d-k}^{-1}$ and  where the distributions on the right-hand sides must be understood as equivalence classes in their respective distribution spaces.
    
    Furthermore, if
    \begin{equation}
        \delta(\dummy - \mat{A}_0) \, r(\dummy - \vec{t}_0) \in \X' = (\X_k \oplus \X_k^\comp)' \subset \Sch'(\Xi_k), 
        \label{eq:filtered-k-plane-ridge-concrete}
    \end{equation}
    where $\X_k$ and $\X_k^\comp$ are complementary Banach spaces and $(\X_k, \X_k')$ form a dual pair of $k$-plane-compatible Banach spaces as in \cref{sec:k-plane-compatible}, then the filtered $k$-plane transform has the \emph{concrete interpretation}
    \begin{equation}
        \KOp_{d-k}\RadonOp_k\curly*{r(\mat{A}_0(\dummy) - \vec{t}_0)} = \P_k^* \curly{\delta(\dummy - \mat{A}_0) \, r(\dummy - \vec{t}_0)} \in \X_k'.
    \end{equation}
\end{theorem}
\begin{proof}
    It suffices to prove the results when $\vec{t}_0 = \vec{0}$ since the general result follows by replacing $r$ with $r(\dummy - \vec{t}_0)$. For all $\varphi \in \Sch(\R^d)$ we have that
    \begin{align*}
        \ang{\delta(\dummy - \mat{A}_0) r(\dummy), \RadonOp_k\curly{\varphi}}_k 
        &= \ang{r, \RadonOp_k\curly{\varphi}(\mat{A}_0, \dummy)} \\
        &= \ang{r(\mat{A}_0(\dummy)), \varphi} \\
        &= \ang{r(\mat{A}_0(\dummy)), \RadonOp_k^*\KOp_{d-k}\RadonOp_k\curly{\varphi}} \\
        &= \ang{\KOp_{d-k} \RadonOp_k r(\mat{A}_0(\dummy)), \RadonOp_k\curly{\varphi}}, \numberthis
    \end{align*}
    where, in the second line, we used \cref{defn:k-plane-ridge-distribution} and, in the third line, we used \cref{thm:cont-inv}. Since $\RadonOp_k\curly{\varphi} \in \Sch_k$, we see from \cref{item:distributional-filtered-k-plane} in \cref{defn:distributional-k-plane} that $\delta(\dummy - \mat{A}_0) r(\dummy) = \KOp_{d-k} \RadonOp_k r(\mat{A}_0(\dummy))$, where the equality is understood in $\Sch_k'$, proving \cref{eq:filtered-k-plane-ridge-abstract}. This then yields \cref{eq:k-plane-ridge-abstract} by substituting $(q * r)$ for $r$. The identity in \cref{eq:filtered-k-plane-ridge-concrete} follows from the application of \cref{thm:k-plane-systematic}.
\end{proof}

\begin{example} \label[example]{example:filtered-k-plane-measure}
    If $r \in \M(\R^{d-k})$, then $\delta(\dummy - \mat{A}_0) \, r(\dummy - \vec{t}_0) \in \M(\Xi_k)$. By \cref{thm:k-plane-compatible}, we know that $\M(\Xi_k) = \M_\iso(\Xi_k) \oplus \paren{\M_\iso(\Xi_k)}^\comp$, where $(C_{0, \iso}(\Xi_k), \M_\iso(\Xi_k))$ forms a dual pair of $k$-plane-compatible Banach spaces. Therefore, by \cref{thm:k-plane-ridges},
    \begin{equation}
        \KOp_{d-k}\RadonOp_k\curly*{r(\mat{A}_0(\dummy) - \vec{t}_0)} = \P_\iso \curly{\delta(\dummy - \mat{A}_0) \, r(\dummy - \vec{t}_0)} \in \M_\iso(\Xi_k). 
    \end{equation}
\end{example}
\begin{remark}
    Let 
    \begin{align*}
        \delta_\iso((\mat{A}, \vec{t}) - (\mat{A}_0, \vec{t}_0)) 
        &\coloneqq \P_\iso\curly{\delta(\dummy - (\mat{A}_0, \vec{t}_0))}(\mat{A}, \vec{t}) \\
        &= \avint_{\mathrm{O}_{d-k}(\R)} \delta((\mat{U}\mat{A}, \mat{U}\vec{t}) - (\mat{A}_0, \vec{t}_0)) \dd\sigma(\mat{U}) \numberthis \label{eq:isotropic-Dirac}
    \end{align*}
    denote the isotropic projection of the Dirac impulse on $\Xi_k$. \Cref{example:filtered-k-plane-measure} establishes the identity
    \begin{equation}
        \KOp_{d-k}\RadonOp_k\curly{\delta(\mat{A}_0(\dummy) - \vec{t}_0)} = \delta_\iso(\dummy - (\mat{A}_0, \vec{t}_0)) \\
        \label{eq:filtered-k-plane-Dirac-ridge}
    \end{equation}
    which is equivalent to
    \begin{equation}
        \RadonOp_k^*\curly{\delta_\iso(\dummy - (\mat{A}_0, \vec{t}_0))} = \delta(\mat{A}_0(\dummy) - \vec{t}_0),
    \end{equation}
    where we use the consequence from \cref{thm:cont-inv-dual} that $\RadonOp_k^* \KOp_{d-k} \RadonOp_k =  \Id$ on $\Sch'(\R^d)$. Furthermore, if $\LOp_d$ is an isotropic LSI operator, then, by the intertwining properties (\cref{thm:intertwine}), we have that
    \begin{equation}
        \LOp_d \RadonOp_k^* \curly{\delta_\iso(\dummy - (\mat{A}_0, \vec{t}_0))} = \RadonOp_k^*\curly{\LOp_{d-k} \delta_\iso(\dummy - (\mat{A}_0, \vec{t}_0))} = r(\mat{A}_0(\dummy) - \vec{t}_0),
        \label{eq:intertwine-k-plane-ridge}
    \end{equation}
    where $r = \FourierOp_{d-k}^{-1}\curly{\hat{L}_{d-k}}$ is the impulse response of $\LOp_{d-k}$.
\end{remark}

%% file: sections/inverse-problems.tex
\section{Application to the Regularization of Inverse Problems} \label{sec:inverse-problems}
In this section, we demonstrate the use of our formalism to specify a new regularization for linear inverse problems and investigate its properties. The resolution of a linear inverse problem involves the estimation of an unknown object $f: \R^d \to \R$ 
based on a finite number of measurements of the form
\begin{equation}
    \vec{y} = \HOp\curly{f} + \vec{\varepsilon} \in \R^M.
\end{equation}
Here, $\HOp$ is a \emph{linear} operator that models the response of the acquisition device (e.g., samples of the Fourier transform in the case of magnetic resonance imaging) and $\vec{\varepsilon} \in \R^M$ is an unknown perturbation or noise term, often assumed to be a vector of i.i.d.\ Gaussian random variables. One typically solves the inverse problem by minimizing a functional of the form
\begin{equation}
    \min_{f \in \X} \: \norm{\vec{y} - \HOp\curly{f}}_2^2 + \lambda \Phi(f),
\end{equation}
where the first term is the \emph{data-fidelity} term and the second term is the \emph{regularization} term that injects prior knowledge on the solution~\cite{EnglRegularization,BenningRegularization}. The parameter $\lambda > 0$ controls the trade off between data fidelity and regularization. The \emph{native space} $\X$ can be thought of as the largest space of functions mapping $\R^d \to \R$ such that the regularizer $\Phi(\dummy)$ is well-defined.

Here, we shall introduce a new sparsity-promoting regularizer in the $k$-plane domain and show that the solutions favor sparse superpositions of $k$-plane ridges. This is a generalization of recently studied regularizers in the Radon domain ($k = (d-1)$) which have been shown to favor neural-network-like solutions~\cite{ParhiShallowRepresenter}. The sparse superpositions of $k$-plane ridges are, in essence, a hybridization of radial basis functions (RBFs) and neural networks.

Given $s \geq 0$, let
\begin{equation}
    \BOp_d^s \coloneqq (\Id - \Delta)^{s/2}
\end{equation}
denote the Bessel potential of order $(-s)$ on $\R^d$, where $\Id$ denotes the identity operator and $\Delta$ denotes the Laplacian. The operator $\BOp_d^s$ can be understood in the Fourier domain via the frequency response
\begin{equation}
    \hat{B}_d^s(\vec{\xi}) = (1 + \norm{\vec{\xi}}_2^2)^{s/2}, \quad \vec{\xi} \in \R^d.
\end{equation}
These operators play an important role in the theory of Sobolev spaces. Indeed, when $s$ is an integer and $1 < p < \infty$, $f \mapsto \norm{\BOp_d^s f}_{L^p}$ is a norm for the Sobolev space $W^{s, p}(\R^d)$. Consider the space
\begin{equation}
    \M_{\RadonOp_k}^s(\R^d) = \curly{f \in \Sch'(\R^d) \st \norm{\KOp_{d-k} \RadonOp_k \BOp_d^s f}_{\M(\Xi_k)} < \infty}.
\end{equation}
When equipped with the norm $\norm{f}_{\M_{\RadonOp_k}^s} \coloneqq \norm{\KOp_{d-k} \RadonOp_k \BOp_d^s f}_{\M(\Xi_k)}$, this is a Banach space. Indeed, the operator
\begin{equation}
    \KOp_{d-k} \RadonOp_k \BOp_d^s: \M_{\RadonOp_k}^s(\R^d) \to \M_\iso(\Xi_k)
\end{equation}
is an isometric isomorphism with inverse given by
\begin{equation}
    \BOp_d^{-s}\RadonOp_k^*: \M_\iso(\Xi_k) \to \M_{\RadonOp_k}^s(\R^d).
    \label{eq:inverse-extreme}
\end{equation}
Since $\M_\iso(\Xi_k)$ is a Banach space, $\M_{\RadonOp_k}^s(\R^d)$ is a Banach space too, when equipped with the previously specified norm. Moreover, it has a ``canonical'' predual which is the space
\begin{equation}
    C_{0, \RadonOp_k}^s(\R^d) = \paren{\KOp_{d-k} \RadonOp_k \BOp_d^s}^* \paren*{C_{0, \iso}(\Xi_k)} = \BOp_d^s \RadonOp_k^* \KOp_{d-k} \paren*{C_{0, \iso}(\Xi_k)},
\end{equation}
i.e., $\paren{C_{0, \RadonOp_k}^s(\R^d)}' = \M_{\RadonOp_k}^s(\R^d)$. This allows us to equip $\M_{\RadonOp_k}^s(\R^d)$ equip with a weak$^*$-topology. With this formalism in hand, we have the following \emph{representer theorem} for $k$-plane-domain regularization.

\begin{theorem} \label{thm:rep-thm}
Let $\HOp\curly{f} = (\ang{h_1, f}, \ldots, \ang{h_M, f}) \in \R^M$ denote the linear measurement functionals where, $h_m \in C_{0, \RadonOp_k}^s(\R^d)$, $m = 1, \ldots, M$, so that the measurement operator is weak$^*$-continuous on $\M_{\RadonOp_k}^s(\R^d)$. Then, for any fixed $\vec{y} \in \R^M$, the solution set to the variational problem
\begin{equation}
    \V \coloneqq \argmin_{f \in \M_{\RadonOp_k}^s(\R^d)} \: \norm{\vec{y} - \HOp\curly{f}}_2^2 + \lambda \norm{\KOp_{d-k} \RadonOp_k \BOp_d^s f}_{\M(\Xi_k)}, \quad \lambda > 0,
\end{equation}
is nonempty, convex, and weak$^*$-compact. The full solution set $\V$ is the weak$^*$ closure of the convex hull of its extreme points, which can all be expressed as
\begin{equation}
    f_\mathrm{extreme}(\vec{x}) = \sum_{n=1}^N a_n \, \rho_s(\mat{A}_n \vec{x} - \vec{t}_n),
\end{equation}
with $N$ parameters (number of atoms with $N \leq M$), $\curly{a_n}_{n=1}^N \subset \R \setminus \curly{0}$, and $\curly{(\mat{A}_n, \vec{t}_n)}_{n=1}^N \subset \Xi_k$ which are data-dependent and not known \emph{a priori}. The profile function $\rho_s: \R^{d-k} \to \R$ is an RBF that corresponds to the Green's function of $\BOp_{d-k}^s$. Finally, the regularization cost, which is common to all solutions, is $\norm{\KOp_{d-k} \RadonOp_k \BOp_d^s f_\mathrm{extreme}}_\M = \sum_{n=1}^N \abs{a_n} = \norm{\vec{a}}_1$.
\end{theorem}
The key takeaway from \cref{thm:rep-thm} is that the solution set to the variational problem is completely characterized by sparse ($N \leq M$) superpositions of \emph{data-adaptive} $k$-plane ridges where the profiles are RBFs with exponential decay~(see \cite[p.~7]{GrafakosModernFourier} for a formula and decay properties for $\rho_s$). By invoking the intertwining properties (\cref{thm:intertwine}), we may rewrite the regularizer as
\begin{equation}
    \norm{\BOp_{d-k}^s \KOp_{d-k} \RadonOp_k f}_{\M(\Xi_k)},
\end{equation}
where we recall that $\BOp_{d-k}^s$ is the $(d-k)$-dimensional Bessel potential of order $(-s)$.

We prove \cref{thm:rep-thm} using recent results characterizing the solution sets to optimization problems over Banach spaces with sparsity-promoting norms~\cite{BoyerRepresenter,BrediesSparsity,UnserUnifyingRepresenter}. These results hinge on the characterization of the extreme points of the unit ball of the regularizer. Since $\M_{\RadonOp_k}^s(\R^d)$ is isometrically isomorphic to $\M_\iso(\Xi_k)$, one can find the extreme points of the unit ball in $\M_{\RadonOp_k}^s(\R^d)$ by applying \cref{eq:inverse-extreme} to the extreme points of the unit ball in $\M_\iso(\Xi_k)$. We characterize these extreme points in \cref{sec:extreme}

\subsection{Extreme Points of the Unit Ball in \texorpdfstring{$\M_\iso(\Xi_k)$}{Miso(Xik)}} \label{sec:extreme}
The canonical evaluation functionals that sample a continuous function defined on $\Xi_k$ at $(\mat{A}_0, \vec{t}_0)$ are the shifted Dirac impulses $\delta(\dummy - (\mat{A}_0, \vec{t}_0)) = \delta(\dummy - \mat{A}_0)\delta(\dummy - \vec{t}_0) \in \M(\Xi_k)$. It is well-known that these correspond to the extreme points of the unit ball in $\M(\Xi_k)$.

We now show that the $k$-plane-compatible sampling functionals are precisely the isotropic projections of the shifted Dirac impulses specified in \cref{eq:isotropic-Dirac} and that these are the extreme points of the unit ball in $\M_\iso(\Xi_k)$.
We summarize their properties in \cref{lemma:Dirac-iso}.
\begin{lemma} \label[lemma]{lemma:Dirac-iso}
The isotropic shifted Dirac impulse $\delta_\iso(\dummy - (\mat{A}_0, \vec{t}_0)) \in \M_\iso(\Xi_k)$ specified by \cref{eq:isotropic-Dirac} satisfies the following properties.
\begin{enumerate}
    \item Sampling: For any $\psi \in C_{0, \iso}(\Xi_k)$,
    \begin{equation}
        \ang{\delta_\iso(\dummy - (\mat{A}_0, \vec{t}_0)), \psi}_k = \psi(\mat{A}_0, \vec{t}_0).
    \end{equation}
    \label{item:k-plane-Dirac-sampling}
    \item Rotation invariance: For any $\mat{U} \in \rmO_{d-k}(\R)$,
    \begin{equation}
        \delta_\iso(\dummy - (\mat{A}_0, \vec{t}_0)) = \delta_\iso(\dummy - (\mat{U}\mat{A}_0, \mat{U}\vec{t}_0)).
    \end{equation}
    \label{item:iso-Dirac-iso}
    \item Unit norm: $\norm{\delta_\iso(\dummy - (\mat{A}_0, \vec{t}_0))}_\M = 1$. \label{item:iso-Dirac-unit-norm}

    \item Linear combination: For any set $\curly{(\mat{A}_n, \vec{t}_n)}_{n=1}^N \subset \Xi_k$ of distinct points,
    \begin{equation}
        \norm*{\sum_{n=1}^N a_n \delta_\iso(\dummy - (\mat{A}_n, \vec{t}_n))}_\M = \sum_{n=1}^N \abs{a_n} = \norm{\vec{a}}_1.
        \label{eq:distinct-ell1}
    \end{equation}
    \label{item:distinct-ell1}
    
    \item Extreme points of  $B_{\M_\iso} = \curly{e \in \M_\iso(\Xi_k) \st \norm{e}_{\M} \leq 1}$: If $e \in \Ext B_{\M_\iso}$, then $e = \pm \delta_\iso(\dummy - (\mat{A}_0, \vec{t}_0))$ for some $(\mat{A}_0, \vec{t}_0) \in \Xi_k$.

    \label{item:Miso-extreme}
\end{enumerate}
\end{lemma}
\begin{proof}
\Cref{item:k-plane-Dirac-sampling,item:iso-Dirac-iso,item:iso-Dirac-unit-norm} are direct consequences of the definition of the isotropic Dirac impulse. \Cref{item:distinct-ell1} follows from the observation that the distributions $\curly{\delta_\iso(\dummy - (\mat{A}_n, \vec{t}_n))}_{n=1}^N$ have disjoint support whenever $\curly{(\mat{A}_n, \vec{t}_n)}_{n=1}^N \subset \Xi_k$ are distinct since this implies equality in the triangle inequality, which proves \cref{eq:distinct-ell1}.

As for \cref{item:Miso-extreme}, the result follows from a classical result~\cite[p. 441]{DunfordLinearOperators} which characterizes the extreme points of the dual of any closed subspace of $\X \subset C_0(Q)$ as the set of all evaluation functionals in $\X'$ multiplied by $\pm 1$, where $Q$ is any (locally compact) Hausdorff space. In the present scenario, we first note that $C_{\iso, 0}(\Xi_k) \subset C_0(\Xi_k)$ is a closed subspace with the domain $\Xi_k$ being a locally compact Hausdorff space. From \cref{item:k-plane-Dirac-sampling}, we know that the evaluation functionals in $\M_\iso(\Xi_k) = \paren*{C_{\iso, 0}(\Xi_k)}'$ are exactly the shifted isotropic Dirac impulses. We refer the reader to~\cite[p. 441]{DunfordLinearOperators} for a proof of this classical result when the domain is a Hausdorff space and to~\cite[Proposition~4.1]{BrediesSparsity} for a proof that applies to any locally compact Hausdorff space.

\end{proof}

\subsection{Proof of \Cref{thm:rep-thm}}
\begin{proof}
    From \cref{item:Miso-extreme} in \cref{lemma:Dirac-iso}, we see that the extreme points of the unit ball
    \begin{equation}
        B_{\M_{\RadonOp_k}^s} \coloneqq \curly{f \in \M_{\RadonOp_k}^s(\R^d) \st \norm{\KOp_{d-k} \RadonOp_k \BOp_d^s f}_\M \leq 1}
    \end{equation}
    take the form $\pm \BOp_d^{-s}\RadonOp_k^* \curly{\delta_\iso(\dummy - (\mat{A}_0, \vec{t}_0))} = \pm \RadonOp_k^* \BOp_{d-k}^{-s} \curly{\delta_\iso(\dummy - (\mat{A}_0, \vec{t}_0))}$, where $(\mat{A}_0, \vec{t}_0) \in \Xi_k$.
    Applying \cref{eq:intertwine-k-plane-ridge} we find that $e \in \Ext B_{\M_{\RadonOp_k}^s}$ if and only if
    \begin{equation}
        e(\vec{x}) = \rho_s(\mat{A}_0\vec{x} - \vec{t}_0), \quad (\mat{A}_0, \vec{t}_0) \in \Xi_k,
    \end{equation}
    where $\rho_s:\R^{d-k} \to \R$ is the impulse response of $\BOp_{d-k}^{-s}$ which is also the Green's function of $\BOp_{d-k}^s$. The abstract representer theorem~\cite[Theorems~2~and~3]{UnserUnifyingRepresenter} then tells us that the solution set $\V$ is nonempty, convex, and weak$^*$-compact and that the full solution set is the weak$^*$ closure of the convex hull of its extreme points. In general, the latter are linear combinations of $N \leq M$ extreme points of the unit regularization ball. In the present scenario, this translates into
    \begin{equation}
        f_\mathrm{extreme}(\vec{x}) = \sum_{n=1}^N a_n \, \rho_s(\mat{A}_n \vec{x} - \vec{t}_n),
    \end{equation}
    where $\curly{a_n}_{n=1}^N \subset \R \setminus \curly{0}$, $\curly{\mat{A}_n}_{n=1}^N \subset V_{d-k}(\R^d)$, $\curly{\vec{t}_n}_{n=1}^N \subset \R^{d-k}$, and $N \leq M$. Moreover, every solution in $\V$ shares the same regularization cost. A calculation in the Fourier domain shows that $\BOp_d^s \rho_s(\mat{A}(\dummy) - \vec{t}) = \delta(\mat{A}(\dummy) - \vec{t})$. Thus, applying the identity \cref{eq:filtered-k-plane-Dirac-ridge}, we find that
    \begin{equation}
         \KOp_{d-k} \RadonOp_k \BOp_d^s f_\mathrm{extreme} = \sum_{n=1}^N a_n \delta_\iso(\dummy - (\mat{A}_n, \vec{t}_n)),
    \end{equation}
    which yields $\norm{\KOp_{d-k} \RadonOp_k \BOp_d^s f_\mathrm{extreme}}_\M = \sum_{n=1}^N \abs{a_n} = \norm{\vec{a}}_1$ (see \cref{item:distinct-ell1}, \cref{lemma:Dirac-iso}).
\end{proof}

%% file: appendix/iso-projector.tex
\section{The Isotropic Projector} \label{app:iso-projector}
In order to derive the explicit formula of the operator that projects a function onto its isotropic part, we first consider projecting functions in $L^2(\R^n)$. This allows us to define the projector as the solution to an optimization problem. Let $L^2_\iso(\R^n)
\subset L^2(\R^n)$ denote the closed subspace of isotropic functions in
$L^2(\R^n)$, as in,
\begin{equation}
  L^2_\iso(\R^n) = \curly{f \in L^2(\R^n) \st f(\mat{U}\vec{x}) = f(\vec{x}),
  \vec{x} \in \R^n, \mat{U} \in \mathrm{O}_n(\R)}.
\end{equation}
The isotropic projection of $f \in L^2(\R^d)$ is then given by its orthoprojection onto $L^2_\iso(\R^n)$ as
\begin{equation}
  \P_\iso f = \argmin_{g \in L^2_\iso(\R^n)} \: \norm{f - g}_{L^2}^2,
  \label{eq:Piso-opt}
\end{equation}
which exists and is unique by the Hilbert projection theorem. We claim that for any $f \in L^2(\R^d)$,
\begin{equation}
  \P_\iso\curly{f}(\vec{x}) = \avint_{\mathrm{O}_n(\R)}
  f(\mat{U}\vec{x}) \dd\sigma(\mat{U}) \in L^2_\iso(\R^n),
\end{equation}
where $\avint$ denotes the average integral, i.e.,
\begin{equation}
  \avint_{\mathrm{O}_n(\R)} \coloneqq \frac{1}{\sigma(\mathrm{O}_n(\R))}
  \,\int_{\mathrm{O}_n(\R)}.
\end{equation}
Clearly, $\P_\iso \P_\iso = \P_\iso$ on $L^2(\R^n)$. Next, for any $f, g \in L^2(\R^n)$,
\begin{align*}
    \ang{\P_\iso f, g}
    &= \int_{\R^n} \paren*{\avint_{\mathrm{O}_n(\R)} f(\mat{U}\vec{x}) \dd \sigma(\mat{U})} g(\vec{x}) \dd\vec{x} \\
    &= \int_{\R^n} \avint_{\mathrm{O}_n(\R)} f(\mat{U}\vec{x}) g(\vec{x}) \dd \sigma(\mat{U}) \dd\vec{x} \\
    &= \avint_{\mathrm{O}_n(\R)} \int_{\R^n} f(\mat{U}\vec{x}) g(\vec{x}) \dd\vec{x} \dd \sigma(\mat{U}) \\
    &= \avint_{\mathrm{O}_n(\R)} \int_{\R^n} f(\vec{y}) g(\mat{U}^\T\vec{y}) \dd\vec{y} \dd \sigma(\mat{U}) \\
    &= \int_{\R^n} \avint_{\mathrm{O}_n(\R)} f(\vec{y}) g(\mat{U}^\T\vec{y}) \dd \sigma(\mat{U}) \dd\vec{y} \\
    &= \int_{\R^n} f(\vec{y}) \paren*{\avint_{\mathrm{O}_n(\R)}  g(\mat{U}^\T\vec{y}) \dd \sigma(\mat{U})} \dd\vec{y} = \ang{f, \P_\iso g}, \numberthis
\end{align*}
where the change of the order of the integrals is justified by the Fubini--Tonelli theorem since
\begin{align*}
    \avint_{\mathrm{O}_n(\R)} \int_{\R^n} \abs{f(\mat{U}\vec{x}) g(\vec{x})} \dd\vec{x} \dd \sigma(\mat{U})
    &\leq \avint_{\mathrm{O}_n(\R)} \norm{f(\mat{U}(\dummy))}_{L^2} \norm{g}_{L^2} \dd\sigma(\mat{U}) \\ &= \norm{f}_{L^2} \norm{g}_{L^2}
    < \infty. \numberthis
\end{align*}
Thus, for all $f \in L^2(\R^n)$, we have that $\P_\iso f$ satisfies the orthogonality principle
\begin{align*}
\ang{f - \P_\iso f, \P_\iso f} 
&= \ang{\P_\iso f - \P_\iso\P_\iso f, f} \\
&= \ang{\P_\iso f - \P_\iso f, f} \\
&= \ang{0, f} = 0. \numberthis
\end{align*}
Therefore, $\P_\iso f$ is indeed the solution to \cref{eq:Piso-opt}.